\documentclass{svjour3}         
\smartqed  
\usepackage{graphicx,rotating}
\usepackage{amsmath,amssymb,enumerate,natbib,color,bbm,ifthen,overpic,ifthen,longtable,lscape}
\usepackage[latin1]{inputenc}

\def\rset{\mathbb R}

\def\eqsp{\;}

\newcommand{\chiopt}[1]{\chi_\mu^\star}

\def\bi{\mathbf{i}}
\def\bj{\mathbf{j}}
\def\g{\mathsf{g}}

\newcounter{rmnum}

\newcommand{\esssup}[2]{\mathop{\mathrm{ess sup}}_{#1} \left( #2\right)}
\newcommand{\essinf}[2]{\mathop{\mathrm{ess inf}}_{#1} \left( #2\right)}
\newcommand{\eqdef}{\ensuremath{\stackrel{\mathrm{def}}{=}}}
\newcommand{\ie}{{\em   i.e.} }

\newcommand{\pp}{\text{a.e.} }
\def\eqsp{\;}
\newcommand{\R}{{\mathbb R}}

\newcommand{\pscal}[2]{\left \langle #1, #2 \right \rangle}
\newcommand{\1}{\ensuremath{\mathbbm{1}}}
\def\Id{\mathrm{Id}}

\newcommand{\strata}[2][]{\ifthenelse{\equal{#1}{}}{\mathsf{S}_{#2}}{\mathsf{S}_{#1,#2}}} 
\def\PE{\operatorname{E}} 
\def\PP{\operatorname{P}} 
\newcommand{\CPE}[3][]
{\ifthenelse{\equal{#1}{}}{\operatorname{E}\left[\left. #2 \, \right| #3 \right]}{\operatorname{E}_{#1}\left[\left. #2 \, \right | #3 \right]}}
\newcommand{\CPP}[3][]
{\ifthenelse{\equal{#1}{}}{\operatorname{P}\left[\left. #2 \, \right| #3 \right]}{\operatorname{E}_{#1}\left[\left. #2 \, \right | #3 \right]}}

\def\PVar{\mathrm{Var}}
\newcommand{\CPVar}[3][]
{\ifthenelse{\equal{#1}{}}{\PVar\left[\left. #2 \, \right| #3 \right]}{\PVar_{#1}\left[\left. #2 \, \right | #3 \right]}}

\def\PP{\mathbb P} 
\def\leb{\lambda} 
\newcommand{\lebh}[2]{\leb_{#1}^{#2}} 
\def\Q{\mathcal{Q}} 
\def\Qopt{\mathcal{Q}^\star}
\newcommand{\qopt}[1]{q^\star_{#1}}
\def\marg{f_\mu} 
\def\mmarg{\zeta_\mu} 
\def\smarg{\psi_\mu} 
\def\obj{V} 

\newlength{\noteWidth}
\newtheorem{rem}{Remark}

\newcounter{hypoconbis}
\newcounter{saveconbis}
\newcommand\debutA{\begin{list} {\textbf{A\arabic{hypoconbis}}}{\usecounter{hypoconbis}}\setcounter{hypoconbis}{\value{saveconbis}}}
\newcommand\finA{\end{list}\setcounter{saveconbis}{\value{hypoconbis}}}

\newcounter{hypocom}
\newcounter{savecom}
\newcommand{\debutB}{\begin{list}{\textbf{B\arabic{hypocom}}}{\usecounter{hypocom}}\setcounter{hypocom}{\value{savecom}}}
\newcommand{\finB}{\end{list}\setcounter{savecom}{\value{hypocom}}}

\def\incg#1#2{\includegraphics[width=#1]{#2}}

\clearpage \newpage

\begin{document}

\title{On adaptive stratification}

\author{Pierre Etore   \and
        Gersende Fort  \and
        Benjamin Jourdain \and
        Eric Moulines
}

\authorrunning{Etore \emph{et al.}} 

\institute{ P. Etore \at
  CMAP, \'Ecole Polytechnique,\\
  Route de Saclay, 91128 Palaiseau Cedex \\
  Tel.: +33 (0)1 69 33 45 67\\
  \email{etore@cmap.polytechnique.fr} \and G. Fort and E. Moulines \at
  Institut des Télécoms, Télécom ParisTech, \\
  46 Rue Barrault, 75634 Paris Cedex 13, France \\
  Tel.: +33 (0)1 45 81 77 82  \\
              \email{surname.name@telecom-paristech.fr}           
           \and
           B. Jourdain  \at
              Université Paris-Est, CERMICS, Projet MathFi ENPC-INRIA-UMLV, \\
              6 et 8 avenue Blaise Pascal, 77455 Marne La Vallée, Cedex 2, France \\
              Tel.: +33 (0)1 64 15 35 67
              \email{benjamin.jourdain@enpc.fr}
}

\maketitle

\begin{abstract}
  This paper investigates the use of stratified sampling as a variance
  reduction technique for approximating integrals over large dimensional
  spaces.  The accuracy of this method critically depends on the choice of the
  space partition, the \emph{strata}, which should be ideally fitted to the
  subsets where the functions to integrate is nearly constant, and on the
  allocation of the number of samples within each strata. When the dimension is
  large and the function to integrate is complex, finding such partitions and
  allocating the sample is a highly non-trivial problem.  In this work, we
  investigate a novel method to improve the efficiency of the estimator "on the
  fly", by jointly sampling and adapting the strata which are hyperrectangles
  and the allocation within the strata. The accuracy of estimators when this
  method is used is examined in detail, in the so-called asymptotic regime
  (\ie\ when both the number of samples and the number of strata are large). It
  turns out that the limiting variance depends on the directions defining the
  hyperrectangles but not on the precise abscissae of their boundaries along
  these directions, which gives a mathematical justification to the common
  choice of equiprobable strata. So, only the directions are adaptively
  modified by our algorithm. We illustrate the use of the method for the
  computation of the price of path-dependent options in models with both
  constant and stochastic volatility. The use of this adaptive technique yields
  variance reduction by factors sometimes larger than 1000 compared to
  classical Monte Carlo estimators.
\end{abstract}

\begin{acknowledgement}
  This work has been written in honour of R.~Rubinstein, for his 70th birthday.
  Most of the ideas used in this paper to reduce the variance of Monte Carlo
  estimator have been inspired by the pioneering work of R.~Rubinstein on
  coupling simulation and stochastic optimization. These very fruitful ideas
  were a constant source of inspiration during this work.

This work is supported by the french National Research Agency (ANR) under the
program ANR-05-BLAN-0299 and by the ``Chair Risques Financiers'', Fondation du Risque.
\end{acknowledgement}

\newpage
\section{Introduction}
A number of problems in statistics, operation research and mathematical finance boils down to the
evaluation of the expectation (or higher order moments) of a random variable
$\phi(Y)$, known to be a complicated real valued function of a vector $Y=
(Y_1,\dots,Y_d)$ of independent random variables. In our applications, we will
mainly focus on simulations driven by a sequence of independent standard normal
random variables, in situations where the dimension $d$ is very large. Such problems arise in
particular in computational finance for the pricing of path-dependent options,
either when the number of underlying assets is large, or when additional source
of randomness is present such as in stochastic volatility models.

The stratification approach consists in dissecting $\rset^d$ into mutually exclusive
\emph{strata} and ensuring that $\phi$ is evaluated for a prescribed
and appropriate number of points in each stratum (see
\cite{glasserman:2004}, \cite{asmussen:glynn:2007}, \cite{rubinstein:kroese:2008}).

The main purpose of this paper is to discuss a way of dissecting the
space into strata and allocating the random draws in the strata,
adapted to the case where $Y$ is a standard Gaussian vector. We also address the accuracy of estimators when this method of sampling
is used, and give conditions upon which the variance reduction is most effective.

Our method makes use of orthogonal directions, to induce a dissection of $\rset^d$ with the right property.
These directions and the associated allocation are learnt adaptively, while the simulations are performed. The advantage of the
adaptive method, similar to those introduced for importance sampling by \cite{rubinstein:kroese:2004} is that information is collected as the simulations are done, and computations of  means and variances of $\phi(Y)$ in strata are used to update  the choice of these strata and of the allocation.
We investigate in some details the asymptotic regime \ie\ where the number of simulations and the number of the strata both go to infinity.

The method is illustrated for pricing path-dependent options driven by high-dimensional Gaussian vectors,
combining importance sampling based on a change of drift together with the suggested
adaptive stratification. The combination of these two methods, already advocated in an earlier work by \cite{glasserman:heidelberger:shahabuddin:1999},
is very effective; nevertheless, these examples show that, contrary to what is suggested in this work, the asymptotical optimal drift vector is not always the most effective direction of stratification.

The paper is organized as follows. In section \ref{sec:stratification}, an introduction to the main ideas of the stratification is presented.
Section~\ref{sec:AsymptoticRegime} addresses the behavior of the stratified estimator in the asymptotic regime (\ie\ when both the number of samples and the number of strata go to infinity).
The roles of the stratification directions, of the strata boundaries in each direction of stratification and of the allocation  within each strata are evidenced.
In section~\ref{sec:AdaptiveAlgorithm}, an algorithm is proposed to adapt the directions of stratification and the allocation of
simulations within each stratum. In Section~\ref{sec:Appli:Finance}, the proposed adaptive stratification procedure is illustrated using applications  for
the pricing of path-dependent options.

\section{An introduction to stratification}
\label{sec:stratification}
Suppose we want to compute an expectation of the form $\PE \left[ \phi(Y)
\right]$ where $\phi : \rset^d \to \rset$ is a measurable function and $Y$ is a
$\rset^d$-valued random variable. We assume hereafter that
\begin{equation}
\label{eq:carint}
\PE\left[ \phi^2(Y) \right] < +\infty \eqsp.
\end{equation}
We consider a \emph{stratification variable} of the form $\mu^T Y$
where $\mu$ is an orthonormal $(d \times m)$ matrix with $m \leq d$. Given a finite partition $ \left\{ \strata{\bi}, \bi \in \mathcal{I} \right\}$ of
$\rset^m$, the sample space $\rset^d$ of $Y$ is
divided into \emph{strata} defined by
\begin{equation}
\label{eq:definition-stratum-mu}
\strata[\mu]{\bi} \eqdef \left\{x \in \rset^d, \mu^T x \in \strata{\bi} \right\}\eqsp, \quad \bi \in \mathcal{I} \eqsp.
\end{equation}
It is assumed in the sequel that the probability of the strata $\{p_\bi, \bi
\in \mathcal{I} \}$
\begin{equation}
  \label{eq:definition-proba-strata}
  p_\bi(\mu) \eqdef \PP\left( Y \in \strata[\mu]{\bi} \right) = \PP\left(\mu^T Y
  \in \strata{\bi} \right) \eqsp,
\end{equation}
are known which is the case when $Y$ is a standard Gaussian vector and
the $S_\bi$ are hyperrectangles. Up to removing some strata, we may assume without loss
of generality that $p_\bi(\mu) > 0$, for any $\bi \in \mathcal{I}$.

Let $M$ be the total number of draws and $\Q =\{q_{\bi}, \bi \in \mathcal{I}\}$
be an allocation vector (i.e. $q_{\bi} \geq 0$ and $\sum_{\bi \in \mathcal{I}}
q_{\bi} =1$)~: the number $M_{\bi}$ of samples allocated to the $\bi$-th
stratum is given by
\begin{equation}
  \label{eq:definition-number-draws}
  M_{\bi} \eqdef \left \lfloor M \sum_{\bj \leq \bi} q_{\bj} \right \rfloor -
\left \lfloor M \sum_{\bj < \bi} q_{\bj} \right \rfloor \eqsp, \quad \bi \in
\mathcal{I} \eqsp,
\end{equation}
where $\lfloor \cdot \rfloor$ denotes the lower integer part and by convention,
$\sum_{\emptyset} q_{\bj} = 0$ (it is assumed that the set of indices
$\mathcal{I}$ is totally ordered). If the number of points in each stratum is chosen to be
proportional to the probability of the strata, the allocation is said to be
\emph{proportional}.
Given the strata $\{ \strata{\bi},
\bi \in \mathcal{I} \}$ and the allocation $\Q$, the \emph{stratified estimator}
with $M$ draws is defined by
\begin{equation}
  \label{eq:definition-stratified-estimator}
  \sum_{\bi \in \mathcal{I}: M_\bi >0}  p_{\bi}(\mu)  \left\{ \frac{1}{M_\bi} \sum_{j=1}^{M_\bi} \phi(Y_{\bi,j})  \right\} \eqsp,
\end{equation}
where $\{Y_{\bi, j}, j \leq M_{\bi}, \bi \in \mathcal{I}\}$ are
independent random variables with $Y_{\bi, j}$ distributed according to the conditional distribution
$\CPP{Y \in \cdot}{\mu^T Y \in
  \strata{\bi}}$ for $j \leq M_{\bi}$.

The stratified estimator is an unbiased estimator of $\PE[\phi(Y)]$ if the $M_{\bi}$'s are all positive (a sufficient condition is $M \geq 1/\min_{\bi} q_\bi $).
Its variance is given by
$\sum_{\bi \in \mathcal{I}: M_\bi >0} \ M_\bi^{-1} p_{\bi}^2(\mu)
\sigma_\bi^2(\mu)$ where $\sigma_{\bi}^2(\mu)$ is the conditional variance of the random vector $\phi(Y)$ given $\mu^T Y
\in \strata{\bi}$,
\begin{equation}
  \label{eq:definition-variance-intrastrata}
  \sigma_\bi^2(\mu) \eqdef \CPE{\phi^2(Y)}{\mu^T Y  \in \strata{\bi}} -  \left(\CPE{\phi(Y)}{\mu^T  Y \in \strata{\bi}} \right)^2 \eqsp.
\end{equation}
When $M$ goes to infinity and the number of strata is either fixed or goes to infinity slowly enough,
the variance of the stratified estimator is equivalent to $M^{-1}
\sum_{\bi \in \mathcal{I}: q_\bi >0} \ q_\bi^{-1} \ p_{\bi}^2(\mu)
\sigma_\bi^2(\mu)$ (see Lemma~\ref{arrondis}). The two key questions that arise in every application of the stratified
sampling method are (i) the choice of the dissection of the space and (ii) for
a fixed $M$, the determination of the number of samples $M_{\bi}$ to be
generated in each stratum $\bi$.  The
optimal allocation minimizing
the above asymptotic variance subject to the constraint $\sum_{\bi \in
  \mathcal{I}} q_{\bi}= 1$ is given by :
\begin{equation}
  \label{eq:optimal-allocation-qi}
  \qopt{\bi}(\mu) \eqdef  \frac{p_\bi(\mu) \; \sigma_\bi(\mu)}{\sum_{\bj \in \mathcal{I}} p_\bj(\mu) \sigma_\bj(\mu)} \eqsp.
\end{equation}

For a given stratification matrix $\mu$, we refer to $\Qopt(\mu) = \{\qopt{\bi}(\mu), \bi \in
\mathcal{I}\}$ as the \emph{optimal stratification vector}.  Of course, contrary to
the proportions $p_\bi(\mu)$, the conditional expectations $\CPE{\phi(Y)}{Y \in
  \strata[\mu]{\bi}}$ are unknown and so are the conditional variances
$\sigma_\bi^2(\mu)$. 

The simplest approach would be to estimate these conditional variances in a pilot run,
to determine the optimal stratification matrix and the optimal allocation vector from these estimates, and then to use them in a second
stage to determine the stratified estimator. Such a procedure is clearly
suboptimal, since the results obtained in the pilot step are not fully
exploited. This calls for a more sophisticated procedure, in the spirit of those
used for adaptive importance sampling; see for example,
\cite{rubinstein:kroese:2004} and \cite{rubinstein:kroese:2008}.  In these
algorithms,  the estimate of conditional variance and the
stratification directions is gradually improved while computing the stratified estimator and estimating its
variance. Such  algorithm extends the procedure by \cite{etore:jourdain:2007},
who proposed to adaptively learn the optimal allocation vector for a set of given strata and
derived a central limit theorem for the adaptive estimator (with the optimal asymptotic variance).

\section{Asymptotic analysis of the stratification performance}
\label{sec:AsymptoticRegime}
We derive in this Section the asymptotic variance of the stratified estimator when both the
total number of draws $M$ and the number of strata (possibly depending upon
$M$) tend to $+\infty$. The variance of the estimator depends on the stratification matrix $\mu$, on the partition $\{ \strata{\bi}, \bi \in \mathcal{I} \}$ of the sample space of $\mu^T Y$ and on the allocation
$\Q$.

For any integer $k$, we denote by $\leb$ the Lebesgue measure on $\rset^k$, equipped with its Borel sigma-field
(the dependence in the dimension $k$ is implicit).  For a probability density $h$
w.r.t the Lebesgue measure on $\rset$, we denote by $H$ its cumulative distribution function,
and $H^{-1}$ its \emph{quantile} function, defined as the generalized inverse
of $H$,
\begin{equation*}
 H^{-1}(u)=  \inf\{x\in \{H>0 \} :H(x)\geq u\} \eqsp, \quad \text{for any $u \in [0,1]$} \eqsp,
\end{equation*}
where, by convention, $\inf\emptyset=+\infty$.  Let $I$ be a positive integer.
The choice of the strata boundaries is parameterized by an $m$-uplet
$(g_1,\dots,g_m)$ of probability densities on $\rset$ in the following sense: for all
$m$-uplet $\bi= (i_1,\dots,i_m) \in \{1, \cdots, I \}^{m}$,
\begin{equation}
\label{eq:definition-stratum}
 \strata{\bi} \eqdef \prod_{k=1}^m \left(G_k^{-1}\left(\frac{i_{k}-1}{I}\right),G_k^{-1}\left( \frac{i_k}{I} \right)\right] \eqsp.
\end{equation}
We denote by $g(x_1,\dots,x_m)\eqdef  \prod_{k=1}^m g_k(x_k)$ the
associated joint density.  Let $\mu$ be a $d \times m$ orthonormal
matrix. We consider the stratification $\mathcal{S}(\mu)= \left\{ \strata[\mu]{\bi}, \bi
  \in \{1,\dots,I\}^m \right\}$ of the space $\rset^d$. Denote by $\varsigma^2_{I,M}(\mu,g,\Q)$ the
  asymptotic variance of the stratified estimator, given by
\begin{equation}
  \label{eq:definition-variance-estimator}
  \varsigma^2_{I,M}(\mu,g,\Q) \eqdef  \sum_{\bi \in \{1,\dots, I\}^m: M_\bi >0}   M_\bi^{-1} \  p_\bi^2(\mu) \sigma_\bi^2(\mu) \eqsp,
\end{equation}
where the number of draws $M_{\bi}$   is given by \eqref{eq:definition-number-draws} and $ p_\bi(\mu),  \sigma^2_\bi(\mu)$, the
probability and the conditional variance are given by (\ref{eq:definition-proba-strata}),  and
(\ref{eq:definition-variance-intrastrata}), respectively. The dependence
w.r.t. $g$ and $\Q$ of $M_\bi$, $p_{\bi}(\mu)$ and $\sigma^2_\bi(\mu)$
is implicit.

We consider allocation vectors $\Q_\chi= \left\{ q_{\bi}(\chi)\eqdef \int_{\strata{\bi}} \chi \ d\leb \eqsp, \bi \in \{1,\dots,I\}^m
\right\}$ parameterized by a probability density $\chi : \rset^m
\to \rset_+$.    We assume that the random
  variable $\mu^T Y$ possesses a density $\marg$ w.r.t.  the Lebesgue measure
  (on $\rset^m$).  We consider the functions
\begin{align*}
  \smarg(x) \eqdef \CPE{\phi^{\phantom{2}}(Y)}{\mu^T Y =x} \eqsp, \quad \text{and} \quad
  \mmarg(x) \eqdef \CPE{\phi^2(Y)}{\mu^T Y =x} \eqsp.
\end{align*}
Using these notations, the asymptotic variance of the stratified estimator may be rewritten as
\[
\varsigma^2_{I,M}(\mu,g,\Q_\chi) = \sum_{\bi \in \{1,\dots, I\}^m: M_\bi >0}
M_\bi^{-1} \ \left\{ \left(\int_{\strata{\bi}}\marg \ d\leb \right)
  \left(\int_{\strata{\bi}}\mmarg\marg \ d\leb
  \right)-\left(\int_{\strata{\bi}}\smarg\marg \ d\leb \right)^2 \right\} \eqsp.
\]
We will investigate the limiting behavior of asymptotic variance $\varsigma^2_{I,M}(\mu,g,\Q_\chi)$ when
the total number of samples $M$ and the number of strata $I$ both tend to
$+\infty$. For that purpose, some technical conditions are required.  For
$\nu$ a measure on $\rset^m$ and $h$ a real-valued measurable function on
$\rset^m$, we denote by $\essinf{\nu}{h}$ and $\esssup{\nu}{h}$ the essential
infimum and supremum w.r.t. the measure $\nu$. From now on we use the following
convention : $z/0$ is equal to $+\infty$ if $z>0$ and to $0$ if $z=0$. \debutA
\item \label{A1} $\int_{\rset^m} \chi^2 /g \, d \leb<+\infty$ and
  $\essinf{g \cdot \leb}{\chi/g} >0$.  \finA
\debutA
\item \label{A2} for $h \in \{ \marg, \mmarg \marg , \smarg \marg \}$,
  $\int_{\rset^m} h^2 /g \ d \leb<+\infty$.  \finA
Under A\ref{A2},
  $\leb$-\pp, $g=0$ implies that $\marg=0$.
Finally, a reinforced integrability condition is needed
\debutA
\item \label{A3} $\int_{\rset^m} \ \marg^4(\mmarg-\smarg^2)^2/[ \chi^2 g] \ d\leb<+\infty$.
\finA
When $m<d$, we establish the expression
of the limit as the number of strata $I$ goes to $+\infty$ of the limiting
variance (as the number of simulations $M$ goes to $+\infty$) of the stratified
estimator. Define
\begin{equation}
  \label{eq:VarLimite}
 \varsigma_\infty^2(\mu,\chi)  \eqdef   \int_{\rset^m} \marg ^2(\mmarg-\smarg^2)/\chi \ d\leb \eqsp.
\end{equation}
\begin{proposition}\label{prop:multVarLimGalCase}
  Let $m$ be an integer such that $m<d$, $g_1, \cdots, g_m$ be probability density functions (pdf) w.r.t. to the Lebesgue measure of $\rset$, $\mu$ be a $d \times m$ orthonormal
  matrix, and  $\chi$ be a pdf w.r.t. the Lebesgue measure on  $\rset^m$.  Assume that $g$ and $\chi$ satisfy assumptions A\ref{A1}-A\ref{A3}.
  Then
\[
\lim_{I \to +\infty}\lim_{M\to +\infty} M\varsigma_{I,M}^2(\mu,g,\Q_\chi) =
\varsigma^2_\infty(\mu,\chi) \eqsp.
\]
Assume in addition one of the following conditions
\begin{enumerate}[(i)]
\item $\esssup{\chi \cdot \leb}{\marg/\chi} <+\infty$ and $\{I_M, M \geq 1
  \}$ is an integer-valued sequence such that $I_M^{-1} + I^m_M M^{-1} \to 0$
  as $M$ goes to infinity.
\item $\{I_M, M \geq 1 \}$ is an integer-valued sequence such that $I_M^{-1} +
  I^{2m}_M M^{-1} \to 0$ as $M$ goes to infinity.
\end{enumerate}
Then,
\[
\lim_{M\to +\infty} M \varsigma_{I_M,M}^2(\mu,g,\Q_\chi) =
\varsigma_\infty^2(\mu,\chi) \eqsp.
\]
\end{proposition}
The proof is given in Section~\ref{sec:Proof:TheoryMultiDim}. It is worthwhile to note that the limiting variance of the stratified estimator
$\varsigma_\infty^2(\mu,\chi)$ does not depend on the densities $(g_1,\dots,g_m)$ that define
the strata : only the
stratification matrix $\mu$ and the allocation vector $\Q_\chi$  enters in the limit.
The contribution to the variance of the
randomness in the directions orthogonal to the rows of $\mu$ dominates at the first order.
In practice, this means that asymptotically, once the directions of
stratification are chosen, the choice of the strata is irrelevant; the usual
choice of $g_i$ as the distribution of the $i$-th component of the
random vector $\mu^T Y$, $i \in \{1,\dots,m\}$ is asymptotically optimal.

On the contrary, the limiting variance $\varsigma_\infty^2(\mu,\chi)$ depends
on the allocation density $\chi$. For a given value of the stratification
directions $\mu$, it is possible to minimize the function $\chi \mapsto
\varsigma_\infty^2(\mu,\chi)$.  Assume that $\int_{\rset^m} \marg
\sqrt{\mmarg-\smarg^2} \ d\leb>0$. Since
$\int_{\rset^m} \marg
\sqrt{\mmarg-\smarg^2} \ d\leb=\PE\left[ \sqrt{\PVar\left[\phi(Y) \vert
        \mu^T Y \right]} \right]\leq \sqrt{\PVar(\phi(Y))},$ the
  integral is finite by \eqref{eq:carint} and it is possible to define a
  density $\chiopt{\mu}$ by
\begin{equation}
\label{eq:definition-chiopt}
\chiopt{\mu} \eqdef \marg \sqrt{\mmarg-\smarg^2}\bigg/\int_{\rset^m} \marg
  \sqrt{\mmarg-\smarg^2} \ d\leb\eqsp.
\end{equation}
Then $\chiopt{\mu}$ is the minimum
of $\chi \mapsto \varsigma_\infty^2(\mu,\chi)$ and the minimal limiting
variance is
\[
\varsigma_\infty^2(\mu,\chiopt{\mu}) = \left(\int_{\rset^m} \marg \sqrt{\mmarg-\smarg^2}
  \ d\leb \right)^2 = \left(\PE\left[ \sqrt{\PVar\left[\phi(Y) \vert
        \mu^T Y \right]} \right] \right)^2\eqsp.
\]
Provided $\chiopt{}$ satisfies assumptions A\ref{A1}-\ref{A2} (note that
in that case, A\ref{A3} is automatically satisfied), the choice
$\chi=\chiopt{}$ for the allocation of the drawings in the strata is
asymptotically optimal.
\begin{rem}
  An expression of the limiting variance $\varsigma_\infty^2(\mu,\chi)$ has been obtained
  in \cite[Lemma 4.1]{glasserman:heidelberger:shahabuddin:1999} in the case $m=1$ and for the proportional allocation
  rule which corresponds to $\chi = \marg$.  It is shown by these authors that the limiting
  variance is $\PE\left(\CPVar{\phi(Y)}{\mu^T Y}\right)$
  which is equal to $\varsigma_\infty^2(\mu, \marg)$ (note that in this case  the stratification density $g=\marg$,
  satisfies the assumptions A\ref{A1}-\ref{A3} provided that $\PE[\phi^4(Y)] < \infty$). Unless $\CPVar{\phi(Y)}{\mu^T Y}$ is a.s.  constant,
  the asymptotic variance is strictly smaller for the optimal choice of the allocation density.
\end{rem}
The optimal allocation density $\chiopt{\mu}$ cannot in general be computed explicitly
but, as shown in the following Proposition, can be approximated by computing the optimal allocation within each stratum.
\begin{proposition}\label{varasopt}
  Let $m < d$ be an integer and $\mu$ be an $(d \times m)$ orthonormal matrix.
  Assume that A\ref{A2} is satisfied.  Then,
\[
\lim_{I \to +\infty} \sum_{\bi \in \{1,\dots, I\}^m} \left| \qopt{\bi}(\mu) -   \int_{\strata{\bi}} \chiopt{\mu} \ d\leb \right| =0 \eqsp,
\]
where $\Qopt (\mu) \eqdef \{\qopt{\bi}(\mu), \bi \in
\{1, \dots, I\}^m \}$ is given by (\ref{eq:optimal-allocation-qi}).  Let
$\{I_M, M \geq 1 \}$ be an integer-valued sequence such that $I_M^{-1} + I^m_M
M^{-1} \to 0$ as $M$ goes to infinity. Then,
\[
\lim_{M \to +\infty} M \varsigma_{I_M, M}^2(\mu,g, \Qopt(\mu)) =
\varsigma_\infty^2(\mu, \chiopt{\mu}) \eqsp.
\]
\end{proposition}
The proof is given in Section~\ref{sec:Proof:TheoryMultiDim}.  As the number of
strata goes to infinity, the stratified estimator run with the optimal
allocation $\Qopt(\mu)$ has the same asymptotic variance as the
stratified estimator run with the allocation deduced from the optimal density
$\chiopt{}$.  In practice, of course, the optimal allocation
$\Qopt(\mu)$ is unknown, but it is possible to construct an
estimator of this quantity by estimating the conditional variance of $\PVar[
\phi(Y) | \mu^T Y \in \strata{\bi}]$ within each stratum \citep {etore:jourdain:2007}.
\begin{remark}When $m=d$, the results obtained are
  markedly different since the accuracy of the stratified estimator now
  depends on the definition of the strata along each direction. Let
$\phi_\mu(x)\eqdef\phi(\mu^T x)$, $\partial_k\phi_\mu$ be the
partial derivative of $\phi_\mu$ w.r.t. its $k$-th coordinate for
$k\in\{1,\dots,d\}$. Let $g_k$ still denote the function
$x=(x_1,\dots,x_d)\in\rset^d\mapsto g_k(x_k)$. Assuming A\ref{A1},
$\esssup{\leb}{f_\mu/g}<+\infty$ and $\phi\in C^1$ satisfies
  $\esssup{\leb}{\sum_{k=1}^d|\partial_k
      \phi_\mu|/g_k}<+\infty$, one checks in \cite{etore:fort:jourdain:moulines:2008} that for any integer-valued sequence $\{I_M, M \geq 1
  \}$ such that  $\lim_{M \to \infty} \left(I_M^{-1} + I^{d+2}_M M^{-1}
  \right) = 0$, $$\lim_{M\to +\infty} MI_M^2 \varsigma_{I_M,M}^2(\mu,g,\Q_\chi) =
\varsigma_\infty^2(\mu,g,\chi)\eqdef\frac{1}{12}\int_{\rset^d}
 \frac{f_\mu^2}{\chi}\sum_{k=1}^d\left(\frac{\partial_k\phi_\mu}{g_k}\right)^2d\leb.$$
In addition, $\lim_{M\to +\infty} MI_M^2\varsigma_{I_M, M}^2(\mu,g, \Qopt(\mu)) =
\varsigma_\infty^2(\mu,g,\chi^\star_{\mu,g})$ with
$\chi^\star_{\mu,g}\propto f_\mu\sqrt{\sum_{k=1}^d\left(\frac{\partial_k\phi_\mu}{g_k}\right)^2}$.
\end{remark}
\section{An adaptive stratification algorithm}
\label{sec:AdaptiveAlgorithm}
As shown in the asymptotic theory presented above, under optimal allocation, it
is more important to optimize the stratification matrix $\mu$ than the strata
boundaries along each stratification direction \footnote{Of course, this is an
  asymptotic result, but our numerical experiments suggest that optimizing the
  strata boundaries along each stratification direction does not lead to a
  significant reduction of the variance. This is why we concentrate on the
  optimization of the stratification matrix}.
Proposition~\ref{prop:multVarLimGalCase} suggests the following strategy: the
``optimal'' matrix $\mu_\star$ is defined as a minimizer of the limiting variance $\mu \mapsto
\varsigma_\infty^2(\mu, \chiopt{})$. Of course, this optimization problem does
not have a closed form expression because the
functions $x \mapsto \smarg(x)$, $x \mapsto \mmarg(x)$ are not available.

We rather use the characterization of $\varsigma_\infty^2(\mu,
\chiopt{})$ as the limiting variance of the
stratified estimator with optimal allocation given in Proposition~\ref{varasopt}. The problem
boils down to search for a minimizer $\mu$ of the variance
$\varsigma^2_{I,M}(\mu,g,\Qopt(\mu))$.  In our applications,  $Y$ is a $d$-dimensional standard normal vector, and $\mu^T Y$ is a $m$-dimensional
standard Gaussian vector. In this case, we set $g_i$, $i=
\{1,\dots,m\}$ to be the standard Gaussian distribution so that the strata boundaries in
each directions are the quantiles of the standard normal variable. Since
$\varsigma_\infty^2(\mu, \chiopt{})$ does not depend on $g$, the impact
of this convenient choice, which leads to equiprobable strata for
the vector $\mu^T Y$, should be limited.

Of course, the optimization of
$\varsigma^2_{I,M}(\mu,g,\Qopt(\mu))$ is a difficult task because
in particular the definition of this function involves multidimensional
integrals, which cannot be computed with high accuracy. Note also that, in most
situations, the optimization should be done in parallel to the main objective,
namely, the estimation of the quantity of interest $\PE[\phi(Y)]$, which is
    obtained using a stratified estimator based on the adaptively defined
    directions of stratification $\mu$. The adaptive stratification is
    analog to the popular adaptive importance sampling; see for example
\cite{rubinstein:kroese:2004},
\cite{arouna:2004}, \cite{kawai:2007}, and \cite{rubinstein:kroese:2008}.

When the function to minimize is an expectation, the classical approaches to
tackle this problem are based on Monte Carlo approximations for
the integrals appearing in the expression of the objective function and its
gradients. There are typically two approaches to Monte Carlo methods, the stochastic approximation procedure and the sample average
approximation method; see \cite{juditsky:lan:nemirovski:shapiro:2007}. In the adaptive stratification context, these Monte Carlo estimators are based on the current fit of the
stratification matrix and of the conditional variances within each stratum, the underlying idea being that
the algorithm is able to progressively learn the optimal stratification, while
the stratified estimator is constructed.

The algorithm described here is closely related to the sample average approximation
method, the main difference with the classical approach being that, at  every time a new search direction
is computed, a new Monte Carlo sample (using the current fit of the strata and of the allocation) is drawn.

Suppose that $Y$ admits a density w.r.t. the Lebesgue measure denoted by $f$.  Define for $\bi
\in \{1, \cdots, I\}^m$, a function $h \in \{f , \phi f , \phi^2 f\}$ ,
and an orthonormal $d \times m$ matrix $\mu$,
\begin{equation}
  \label{eq:FunctionNuI}
  \nu_\bi(h,\mu) \eqdef \int_{\strata[\mu]{\bi}} h   \ d\leb   =  \int \prod_{k=1}^m \1_{\{y, G_k^{-1}((\bi_k-1)/I) \leq \pscal{\mu_k}{y}  \leq G_k^{-1}(\bi_k/I)\}}    h   \ d\leb   \eqsp,
\end{equation}
where $\pscal{x}{y}$ denotes the scalar product of the vectors $x$ and
$y$ and $\mu_k$ the $k$-th column of $\mu$.
Using the definition of $\nu_\bi$, the proportions $p_\bi(\mu)$ and the
conditional variances in each stratum $\sigma_\bi^2(\mu)$ respectively given
by \eqref{eq:definition-proba-strata} and~
\eqref{eq:definition-variance-intrastrata} may be expressed as, when $\nu_{\bi}(f,\mu) > 0$,
\begin{equation}
  \label{eq:PandSigma}
  p_\bi(\mu)=\nu_\bi(f,\mu) \eqsp, \quad \text{and} \quad
\sigma_\bi^2(\mu)=\frac{\nu_\bi(f
  \phi^2,\mu)}{\nu_\bi(f,\mu)}-\left(\frac{\nu_\bi(f
    \phi,\mu)}{\nu_\bi(f,\mu)}\right)^2 \eqsp.
\end{equation}
When $M$ is large and $I$ is fixed, minimizing the asymptotic variance of the
stratified estimate with optimal allocation is equivalent to minimize
$\obj(\mu)$ w.r.t. the stratification matrix $\mu$ where (see
Lemma~\ref{arrondis})
\[
\obj(\mu) \eqdef \sum_{\bi=1}^\mathcal{I} p_\bi(\mu) \sigma_\bi(\mu) =
\sum_{\bi=1}^\mathcal{I} \left( \nu_\bi(f,\mu) \nu_\bi(f \phi^2,\mu) -
  \nu_\bi^2(f \phi, \mu)\right)^{1/2} \eqsp.
\]
Assuming that the functions $\mu \mapsto \nu_\bi(h,\mu)$ are differentiable at $\mu$ for $h \in \{f,f
\phi, f \phi^2\}$, the gradient
may be expressed as
\begin{equation}
  \label{eq:ExpressionGradV}
    \nabla_{\mu} \; V(\mu) = \sum_{\bi=1}^\mathcal{I} \frac{ \nabla_\mu \nu_\bi(f,\mu) \
    \nu_\bi(f\phi^2,\mu) + p_{\bi}(\mu) \ \nabla_{\mu} \nu_\bi(f\phi^2,\mu) - 2
    \nu_\bi(f\phi,\mu) \ \nabla_{\mu} \nu_\bi(f\phi,\mu) }{2 \ p_{\bi}(\mu)
    \sigma_{\bi}(\mu) } \1_{\{p_{\bi}(\mu) \sigma_{\bi}(\mu) \neq 0\}} \eqsp.
\end{equation}
The computation of this gradient thus requires to calculate $\nabla_\mu \, \nu_\bi(h,\mu)$ for $h \in \{f, f \phi, f \phi^2 \}$.
For a vector $\nu \in \rset^d$, $\nu \neq 0$, and $z \in \rset$, define
$\lebh{z}{\nu}$, the restriction of the Lebesgue measure on the hyperplane $\{y \in \rset^d,
\pscal{\nu}{y} = z \}$.
\begin{proposition}
  \label{prop:DeriveeMu}
  Let $z\in\rset$, $h:\rset^d\rightarrow\R$ be a locally bounded
  integrable real function, $g_z:\rset^d\ni\nu\mapsto \int
  \1_{\{y, \pscal{\nu}{y} \leq z \} } h(y) \ d\leb(y) \eqsp$ and $\mu \in \rset^d$ be a non-zero vector. Assume that $h$ is continuous
  $\lebh{z}{\mu}$ almost everywhere and that there exists $\varepsilon>0$ such
  that
  \begin{equation}
  \label{eq:uniforme-integrabilite}
  \lim_{M\rightarrow +\infty}\sup_{ |\nu - \mu| \leq \varepsilon}\int |y| \1_{\{|y|\geq M\}}|h(y)|  \; d\lebh{z}{\nu}(y)=0 \eqsp.
  \end{equation}
  Then, the function $\nu \mapsto g_z(\nu)$ is differentiable at $\mu$ and
  $\nabla_\mu \; g_z(\mu)=- \int \frac{y}{|\mu|} \; h( y) \;d\lebh{z}{\mu}(y)$.\end{proposition}
\begin{corollary}
\label{coro:Derivee-Matrice}
Assume that $h$ is a real locally bounded integrable function. Let $m$
be an integer, $z= (z_1,\dots,z_m) \in \rset^m$, $g_z: \rset^{d \times m}\ni (\nu_1,\dots,\nu_m)\mapsto \int \prod_{k=1}^m \1_{\{y, \pscal{\nu_k}{y} \leq z_k \}} h(y) \ d\leb(y)$
and $\mu= [\mu_1,\dots, \mu_m] \in \rset^{d \times m}$ be a full rank matrix.
Assume that $h$ is continuous $\sum_{k=1}^m\lebh{z}{\mu_k}$ almost
everywhere and that there exists $\varepsilon>0$ such that, for any $k \in \{1,\dots,m\}$,
$\lim_{M\rightarrow +\infty}\sup_{ |\nu - \mu_k| \leq \varepsilon}\int |y| \1_{\{|y|\geq M\}}|h(y)|  \; d\lebh{z}{\nu}(y)=0$.
Then, $g_z$ is differentiable at $\mu$ and the differential $\nabla_\mu g_z$ is given by $\nabla_\mu g_z = [\nabla_{\mu_1} g_z, \dots, \nabla_{\mu_m} g_z]$, where
$$
\nabla_{\mu_i} g_z(\mu) = -\int \frac{y}{|\mu_i|} \;\prod_{k \ne i} \1_{\{y, \pscal{\mu_k}{y} \leq z_k \}}  h( y) \;d\lebh{z_i}{\mu_i}(y) \eqsp.
$$
\end{corollary}

The algorithm goes as follows. Denote by $\{ \gamma_t \}$ a sequence of stepsizes.
Consider the strata $\{\strata{\bi}, \bi \in \{1, \cdots, I \}^m \}$ given by
(\ref{eq:definition-stratum}) for some product density $g$.
\begin{list}{}{}
\item \begin{enumerate}
  \item {\sf Initialization.} Choose initial stratification directions
    $\mu^{(0)}$ and an initial number of draws in each statum $M^{(0)} \eqdef
    \{M_\bi^{(0)}, \bi \in \eqdef \{1, \dots, I \}^m\}$ such that $\sum_\bi
    M_\bi^{(0)} = M$. Compute the probabilities $p_\bi(\mu^{(0)})$ of each
    stratum.
  \item {\sf Iteration.} At iteration $t+1$, given $\mu^{(t)}$, $M^{(t)}$ and
    $\{p_\bi(\mu^{(t)}), \bi \in \{1, \cdots, I\}^m\}$,
\begin{enumerate}
\item \label{algo:step1} {\em Compute $\widehat{\nabla \obj}(\mu^{(t)})$:}
  \begin{enumerate}[(i)]
  \item \label{algo:step1a} for $ \bi \in \{1,\cdots, I\}^m$, draw
    $M_{\bi}^{(t)}$ realizations of i.i.d.  random variables
    $\{Y_{\bi,k}^{(t)}, k \leq M_{\bi}^{(t)}\}$ with distribution $\PP(Y \in
    \cdot \vert Y \in \strata[\mu^{(t)}]{\bi})$ and evaluate $\hat{\nu}_\bi^{(t+1)}(h) = \frac{p_\bi(\mu^{(t)})}{M_{\bi}^{(t)}}
\sum_{k=1}^{M_{\bi}^{(t)}} \ h\left(Y_{\bi,k}^{(t)} \right)$ for $h \in \{
    \phi,  \phi^2 \}$.
\item \label{algo:step1b} for $k \in \{1, \cdots, m\}$, $ s \in
  \{G_k^{-1}(1/I), \cdots, G_k^{-1}((I-1)/I)\}$, draw $\tilde M_{k,s}^{(t)}$
  realizations of i.i.d.  random variables with distribution $\PP(Y \in \cdot
  \vert [\mu^{(t)}_k]^T Y = s)$.  Compute a Monte Carlo estimate of
  $\nabla_{\mu}\nu_\bi(h,\mu^{(t)})$ for $h \in \{f, f\phi, f\phi^2 \}$ based
  on Corollary~\ref{coro:Derivee-Matrice}.
\item Compute a Monte Carlo estimate of
  $\nabla \obj(\mu^{(t)})$ based on the expression (\ref{eq:ExpressionGradV}).
\end{enumerate}
\item {\em Update the direction of stratification:} Set $ \tilde \mu =\mu^{(t)}
  - \gamma_t \ \widehat{\nabla \obj}(\mu^{(t)})$; define $\mu^{(t+1)}$ as the
  orthonormal matrix found by computing the singular value decomposition of
  $\tilde{\mu}$ and keeping the $m$ left singular vectors.
\item {\em Update the allocation policy:}
  \begin{enumerate}[(i)]
  \item compute an estimate $\hat{\sigma}_\bi^{(t+1)}$ of the standard
    deviation within stratum $\bi$
\[
\hat{\sigma}_\bi^{(t+1)} = \left( \frac{\hat{\nu}_\bi^{(t+1)}(\phi^2)}{p_\bi(\mu^{(t)})}-\left(\frac{\hat{\nu}_\bi^{(t+1)}(\phi)}{p_\bi(\mu^{(t)})}\right)^2 \right)^{1/2}\eqsp.
\]
\item  Update the allocation vector
  \[
  q_\bi^{(t+1)} = \frac{p_\bi(\mu^{(t)}) \ \hat{\sigma}_\bi^{(t+1)}}{\sum_{\bj
      \in \{1, \dots, I \}^m } p_\bj(\mu^{(t)}) \ \hat{\sigma}_\bj^{(t+1)}}
  \eqsp,
  \]
  and the number of draws $\{M_\bi^{(t+1)}, \bi \in \{1, \dots, I \}^m \}$ by
  applying the formula (\ref{eq:definition-number-draws}) with a total number
  of draws equal to $M$.
 \end{enumerate}
\item {\em Update the probabilities} $p_\bi(\mu^{(t+1)})$, $\bi \in \{1,\cdots,
  I\}^m$.
\item \label{algo:stepz}{\em Compute an averaged stratified estimate of the
    quantity of interest: } Estimate the Monte Carlo variance of the stratified
  estimator for the current fit of the strata and the optimal allocation
\[
[\varsigma^2]^{(t+1)} = \frac{1}{M} \left( \sum_{\bi \in \{1,\cdots, I\}^m}
  p_\bi(\mu^{(t)}) \ \hat{\sigma}_\bi^{(t+1)} \right)^2 \eqsp.
\]
Compute the current fit of the stratified estimator by the following weighted average
\begin{equation}
  \label{eq:AveragedEstimate}
  \mathcal{E}^{(t+1)} = \left( \sum_{\tau =1}^{t+1}
  \frac{1}{[\varsigma^2]^{(\tau)}} \right)^{-1} \ \sum_{\tau =1}^{t+1}
\frac{1}{[\varsigma^2]^{(\tau)}} \sum_{\bi \in \{1,\cdots, I\}^m}
\hat{\nu}_\bi^{(\tau)}(\phi) \eqsp.
\end{equation}
\end{enumerate}
\end{enumerate}
\end{list}
There are two options to choose the stepsizes $\{\gamma_t, t \geq 0\}$. The
traditional approach consists in taking a decreasing sequence satisfying the
following conditions (see for example \cite{pflug:1996,kushner:yin:2003})
\[
\sum_{t \geq 0} \gamma_t = +\infty \eqsp, \qquad \qquad \sum_{t \geq 0}
\gamma_t^2 < +\infty \eqsp.
\]
If the number of simulations is fixed in advance, say equal to $N$, then one can use a constant stepsize strategy,
\ie\ choose $\gamma_t= \gamma$ for all $t \in \{1, \dots, N\}$. As advocated in \cite{juditsky:lan:nemirovski:shapiro:2007}, a
sensible choice in this setting is to take $\gamma_t$ proportional to $N^{-1/2}$.
The convergence of this crude gradient algorithm proved to be quite fast in
all our applications, so it is not required to resort to computationally
involved alternatives.

Step~\ref{algo:step1b} is specific to the optimization problem to solve and is
not related to the stratification itself. The number of draws for the
computation of the surface integral (see Corollary~\ref{coro:Derivee-Matrice})
can be chosen independently of the allocation $M^{(t)}$.  When the samples in
steps~\ref{algo:step1a} and \ref{algo:step1b} can be obtained by transforming
the same set of variables (see Section~\ref{sec:Appli:Finance} for such a
situation), it is natural to choose $\tilde M^{(t)} = \{\tilde M_{k,s}^{(t)}, k
\in \{1,\cdots, m\}, s \in \{G_k^{-1}(1/I), \cdots, G_k^{-1}((I-1)/I)\}\}$ such
that $\sum_{k,s} \tilde M_{k,s}^{(t)} =M$.

When $\marg$ has a product form (which is the case e.g. when $Y$ is a standard
d-dimensional Gaussian distribution), we can set $g = \marg$. Then, the strata
are equiprobable and $p_\bi(\mu) =1/I^m$ for any ($\bi$, $\mu$).

It is out of the scope of this paper to prove the convergence of this algorithm and we refer the reader to
classical treatises on this subject.
The above algorithm provides, at convergence, both \textit{(i)} ``optimal''
directions of stratification and an estimate of the associated optimal
allocation; \textit{(ii)} an averaged stratified estimate $\mathcal{E}$. By
omitting the step~\ref{algo:stepz}, the algorithm might be seen as a mean for
computing the stratification directions and the associated optimal allocation,
and these quantities can then be plugged in a ``usual'' stratification
procedure.

\section{Applications in Financial Engineering}
\label{sec:Appli:Finance}
The pricing of an option amounts to compute the expectation $ \PE\left[\Xi(Y)
\right]$ for some measurable non-negative function $\Xi$ on $\rset^d$, where
$Y$ is a standard $d$-multivariate Gaussian variable.  The Cameron-Martin
formula implies that for any $\nu\in \rset^d$,
\begin{equation}
  \label{eq:Cameron-Martin}
\PE\left[\Xi(Y) \right]  = \PE\left[ \Xi(Y+\nu) \ \exp( -\nu^T Y - 0.5 \nu^T
  \nu)\right] \eqsp,
\end{equation}
The variance of the plain Monte Carlo estimator depends on the choice of $\nu$.
In all the experiments below (except for the
\cite{glasserman:heidelberger:shahabuddin:1999} estimator), we use either $\phi(y) =
\Xi(y)$ (case $\nu=0$) or $\phi(y) = \Xi(y+\nu_\star) \ \exp( -\nu^T_\star y -
0.5 \nu^T_\star \nu_\star)$ where $\nu_\star$ is the solution of the
optimization problem
\begin{equation}
  \label{eq:Caracterization-Drift-Nu}
  \mathrm{argmax}_{\{\nu \in \rset^d, \Xi(\nu)>0\}} \ \ \left\{ \ln \Xi(\nu) -
  0.5 \nu^T \nu \right\}\eqsp,
\end{equation}
(case $\nu = \nu_\star$).  The motivations for this particular choice of the
drift vector $\nu$ and procedures to solve this optimization problem are
discussed in~\cite{glasserman:heidelberger:shahabuddin:1999}.

We apply the adaptive stratification procedure introduced in
Section~\ref{sec:AdaptiveAlgorithm} (hereafter referred to as ``{\tt
  AdaptStr}'') in the case $m=1$. For comparison purposes, we also run the
stratification procedure proposed
in~\cite{glasserman:heidelberger:shahabuddin:1999} (hereafter referred to as ``
{\tt GHS}''), combining \textit{(i)} importance sampling with the drift
$\nu_\star$ defined in \eqref{eq:Caracterization-Drift-Nu}, and \textit{(ii)}
stratification with proportional allocation and direction $\mu_\g$ defined
in~\cite[Section~4.2]{glasserman:heidelberger:shahabuddin:1999}.

We also run stratification algorithms with three different directions of
stratification: the vector $\mu_\star \propto \nu_\star$, the vector {\tt $
  \mu_{\mathrm{reg}}$} proportional to the vector of linear regression of the function $\phi(Y)$ on
$Y$ (these regression coefficients are obtained in a pilot run), and a vector $\mu_l$ which is
a simple guess specific to each application. For these three directions, we run stratification
with proportional allocation (case ``$q_i$'' set to ``\emph{prop}'') and with
optimal allocation (case ``$q_i$'' set to ``\emph{opt}''). We also run the
plain Monte Carlo estimator (column ``{\tt MC}''); when used with
$\nu=\nu_\star$, ''{\tt MC}'' corresponds to an importance sampling estimator
with a drift function $\nu_\star$.

Finally, we compare these stratified estimators to Latin Hypercube (LH)
estimators (see \cite{glasserman:2004}, \cite{owen:2003} for a description of
this method).  For $Y$ a standard normal vector in $\rset^d$, the expectation
of interest $\PE[\phi(Y)]$ is also equal to $\PE[\phi(OY)]$ for any orthogonal
matrix $O\in\rset^{d\times d}$ but the variance of the LH estimator associated
with the variable $\phi(OY)$  depends on the choice of $O$.  Unfortunately, it is
very difficult to compute explicitly the asymptotic variance of LH estimators
and therefore to adapt the matrix $O$; see \cite{owen:1992}. 
Since LH somehow consists in stratifying each canonical direction,
choosing the first column of $O$ equal to the
stratification direction $\mu$ should be sensible.
In our numerical experiments, we consider such matrices $O$ obtained by orthonormalization of the basis combining $\mu$ and the
  $d-1$ last vectors of the canonical basis of $\R^d$ with $\mu$ equal
  to $\mu_{\star}$, $\mu_{\mathrm{reg}}$ or to the adaptive stratification direction obtained by our algorithm {\tt AdaptStr}.

\subsection{Practical implementations of the  adaptive stratification procedure}
\label{subsubsec:implementation}
The numerical results have been obtained by running {\tt Matlab} codes
available from the authors \footnote{These
codes are freely available from the url {\tt
  http://www.tsi.enst.fr/$\sim$gfort/}}
In the numerical applications below, $m=1$. We choose $g = \marg$ so that the
strata are equiprobable ($p_\bi(\mu) = 1/I$).  We choose $I = 100$ strata and
$M =20\,000$ draws per iterations.

The drift vector $\nu$ that solves (\ref{eq:Caracterization-Drift-Nu}) is
obtained by running {\tt solnp}, a nonlinear optimization program in Matlab
freely available at {\tt http://www.stanford.edu/$\sim$yyye/matlab/}.  The
direction $\mu^{(0)}$ is set to the unitary constant vector $(1, \cdots,
1)/\sqrt{d}$; the initial allocation $M^{(0)}$ is proportional.  Exact sampling
under the conditional distributions $\PP(Y \in \cdot \vert Y \in
\strata[\mu^{(t)}]{\bi})$ and $\PP(Y \in \cdot \vert [\mu^{(t)}]^T Y =s)$ can
be done by linear transformation of standard Gaussian vectors (see
\citep[section 4.3, p. 223]{glasserman:2004}).  The draws in step
\ref{algo:step1a} and \ref{algo:step1b} can be obtained by transforming the
same set of $M^{(t)}$ Gaussian random variables $\{V^\bi_j, j \leq M_\bi^{(t)},
\bi \in \{1, \cdots,I\}\}$.  Therefore, the total number of $d$-dimensional
Gaussian draws by iteration is $M$ (the estimates of $\nu_\bi(h,\mu)$ and
$\nabla_\mu \nu_\bi(h,\mu)$ are not independent); $M$ uniform draws in $(0,1)$
are also required to sample under the conditional distribution $\PP(Y \in \cdot
\vert Y \in \strata[\mu^{(t)}]{\bi})$. The criterion is optimized using a fixed
stepsize steepest descent algorithm (the stepsize is determined using a limited
set of pilot runs).

\subsection{Assessing efficiency of the adaptive stratification procedure}
\label{sec:assessing_efficiency}
We compare the averaged stratified estimate $\mathcal{E}^{(N)}$ obtained after
$N=200$ iterations, with different stratification procedures and with the crude
Monte Carlo estimate.  We report in the tables below the estimate of the option
prices and the estimates of the variance of the estimator obtained from 50
independent replications.

The comparison of the procedures relies on the variance of the estimators.  The
column ``{\tt MC}'' is an estimate of the variance of $\phi(Y)$ computed with
$MN$ i.i.d. samples of a $d$-multivariate gaussian distribution. In the case
$\nu = 0$, this is an estimation of the variance of the plain Monte Carlo
estimator; when $\nu = \nu_\star$, this corresponds to an estimation of the
Importance Sampling estimator (with importance sampling distribution equal to a
standard Gaussian distribution centered at $\nu_\star$). The column ``{\tt
  AdaptStr}'' is the limiting variance per sample of $\mathcal{E}^{(N)}$ which is
equal to
\[
N \, \left\{ \sum_{t=1}^N \left(\left[\sum_{\bi} p_\bi \,
      \hat{\sigma}_\bi^{(t)}\right]^2 \right)^{-1} \right\}^{-1} \sim
\left(\sum_{\bi} p_\bi \, \sigma_\bi(\mu^{(+\infty)}) \right)^2\eqsp,
\]
when each iteration $t \in \{1, \cdots, N\}$ implies $M$ draws (see the
algorithm in Section~\ref{sec:AdaptiveAlgorithm}); note that by definition of
our procedure, the allocation is optimal.  The column ``{\tt GHS}'' is an
estimate of $\sum_\bi p_\bi \sigma_\bi^2(\mu_\g)$ computed with $MN$ samples;
note also that by definition of this procedure, only the case $\nu=\nu_\star$
and the proportional allocation has been considered. The columns ``{\tt
  $\mu_{\mathrm{reg}}$}'', ``{\tt $\mu_{\star}$}'', and ``{\tt $\mu_{l}$}''
report the results for the stratification procedures with these directions of
stratification: the rows 'proportional allocation' report an estimation of
$\sum_\bi p_\bi \sigma_\bi^2(\mu)$ computed with $MN$ samples (for $\mu \in
\{\mu_{\mathrm{reg}}, \mu_\star, \mu_l\}$ and $\nu \in \{0, \nu_\star \}$). We
also consider the results for the optimal allocation, and to that goal we
estimate the standard deviation within each stratum by an iterative algorithm -
with $N$ iterations - : the rows 'optimal allocation' report an estimation of
$N \, \left\{ \sum_{t=1}^N \left([\sum_{\bi} p_\bi \,
    \hat{\sigma}_\bi^{(t)}(\mu)]^2 \right)^{-1} \right\}^{-1}$ where
$\{\hat{\sigma}_\bi^{(t)}(\mu), \bi \leq I\}$ is an estimation of the standard
deviation of the strata computed with a total number of $M$ draws allocated to
each stratum according to the optimal allocation computed at the previous
iteration $(t-1)$ (the allocation at iteration $0$ is the proportional one).
For Latin Hypercube samplers, the total number of draws ($M N$) are allocated
to generate $N$ i.i.d. estimators $\mathcal{E}^{(k)}_M$, $k \in \{1, \cdots,
N\}$, each based on a Latin Hypercube sample of size $M$. The estimate {\tt
  LHS} is the average of these $N$ estimators; we also report the variance
equal to $ M \ \left\{ N^{-1} \sum_{k=1}^N [\mathcal{E}^{(k)}_M]^2 - \{N^{-1}
  \sum_{k=1}^N \mathcal{E}^{(k)}_M\}^2 \right\}$.

\subsection{Asian options}
\label{subsec:asian}
Consider the pricing of an arithmetic Asian option on a single underlying asset
under standard Black-Scholes assumptions. The price of the asset is described
by the stochastic differential equation $\frac{dS_t}{S_t} = r \, dt + \upsilon
\, d W_t \eqsp, \quad S_0= s_0,$ where $\{W_t, t \geq 0\}$ is a standard
Brownian motion, $r$ is the risk-free mean rate of return, $\upsilon$ is the
volatility and $s_0$ is the initial value. The asset price is discretized on a
regular grid $0=t_0< t_1 < \cdots < t_d=T$, with $t_i \eqdef i T/d$. The
increment of the Brownian motion on $[t_{i-1},t_i)$ is simulated as $\sqrt{T/d}
\, Y_i$ for $i \in \{1,\cdots, d\}$ where $Y=(Y_1, \cdots, Y_d) \sim
\mathcal{N}_d(0, \Id)$. The discounted payoff of a discretely monitored
arithmetic average Asian option with strike price $K$ is given by $\Xi(Y)$,
\[
\Xi(y) =\exp(-rT) \left( \frac{s_0}{d} \sum_{k=1}^d \exp\left((r - 0.5
    \upsilon^2) \frac{k T}{d} + \upsilon \sqrt{\frac{T}{d}} \sum_{j=1}^k y_j
  \right) -K \right)_+ \eqsp \eqsp, \quad y=(y_1,\cdots, y_d) \in \rset^d
\eqsp,
\]
where for $x \in \rset$, $x_+ = \max(x,0)$.  In the numerical applications, we
take $s_0 = 50$, $r=0.05$, $T=1$, $(\upsilon,K) \in \{(0.1,45), (0.5,45),
(0.5,65), (1,45), (1,65) \}$ and $d = 16$.  We choose $\mu_l \propto (d,d-1,
\cdots, 1)$.

\begin{figure}
   \centering
   \incg{200pt}{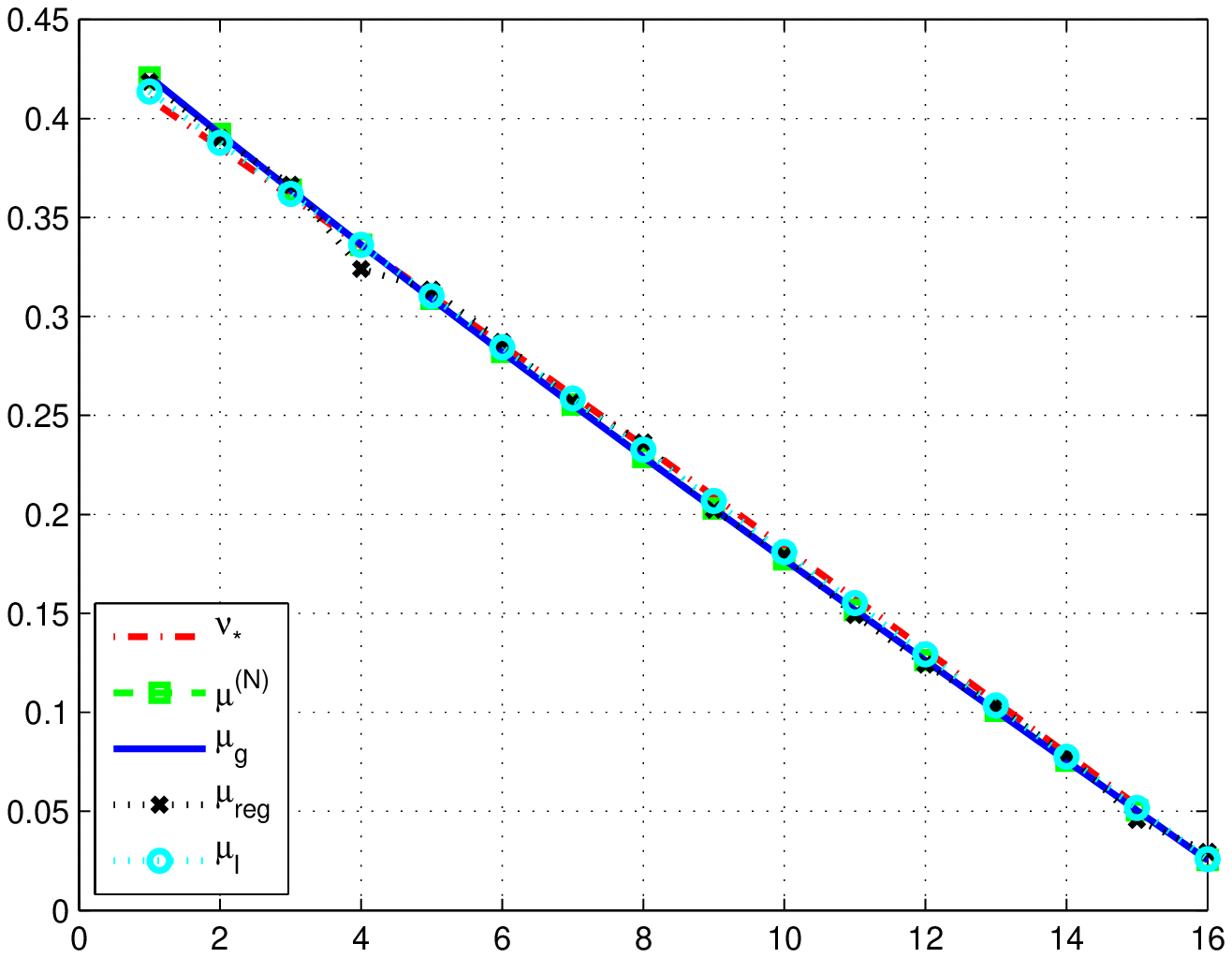}  \incg{200pt}{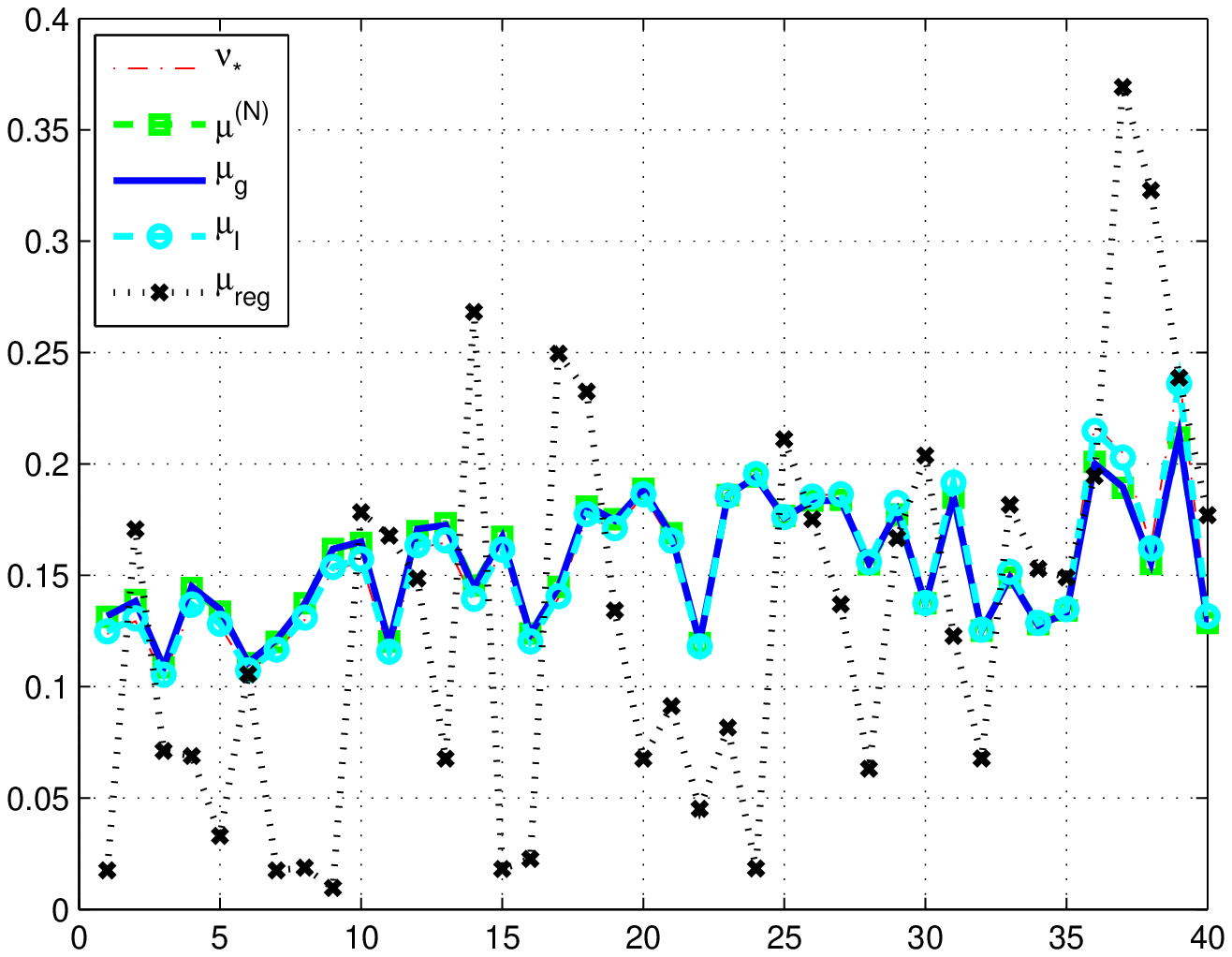}
   \caption{[left panel] Asian Option when $(\upsilon, K,\nu) = (0.1,45,\nu_\star)$:  drift vector $\nu_\star$ and  directions $\mu^{(N)}$, $\mu_\g$, $\mu_{\mathrm{reg}}$ and $\mu_l$.  $\nu_\star$ has been scaled to have norm $1$
     ($\nu_\star \leftarrow \nu_\star / 0.42$). [right panel] Basket Option when
     $(c,K,\nu) = (0.1, 45, \nu_\star)$: drift vector $\nu_\star$ and
     directions $\mu^{(N)}$, $\mu_\g$, $\mu_l$, $\mu_{\mathrm{reg}}$ and
     $\mu_l$.  $\nu_\star$ has been scaled to have norm $1$ ($\nu_\star
     \leftarrow \nu_\star / 0.41$).}
 \label{fig:Asie:directions}
 \end{figure}

\begin{figure}[h]
   \centering \incg{200pt}{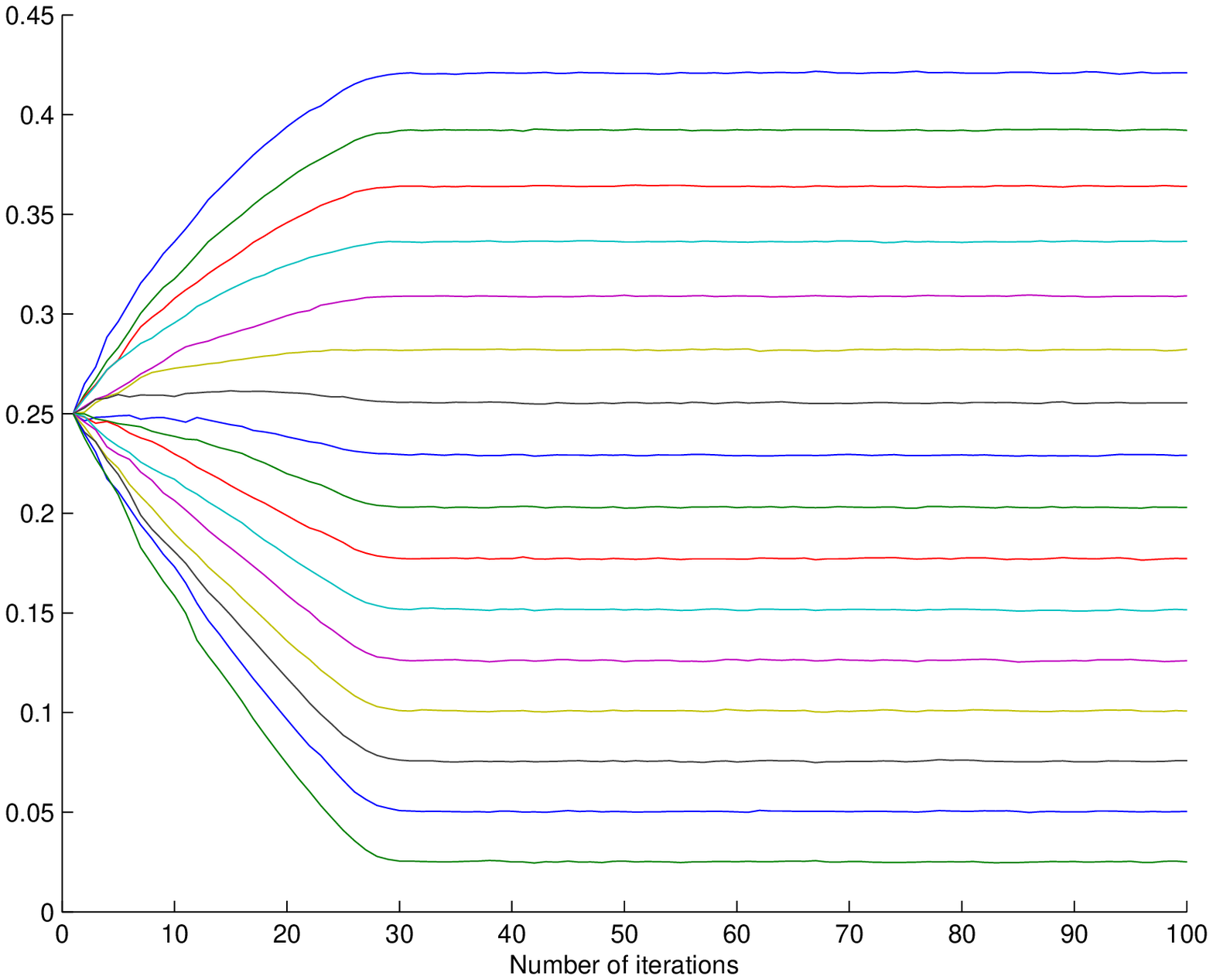}  \incg{200pt}{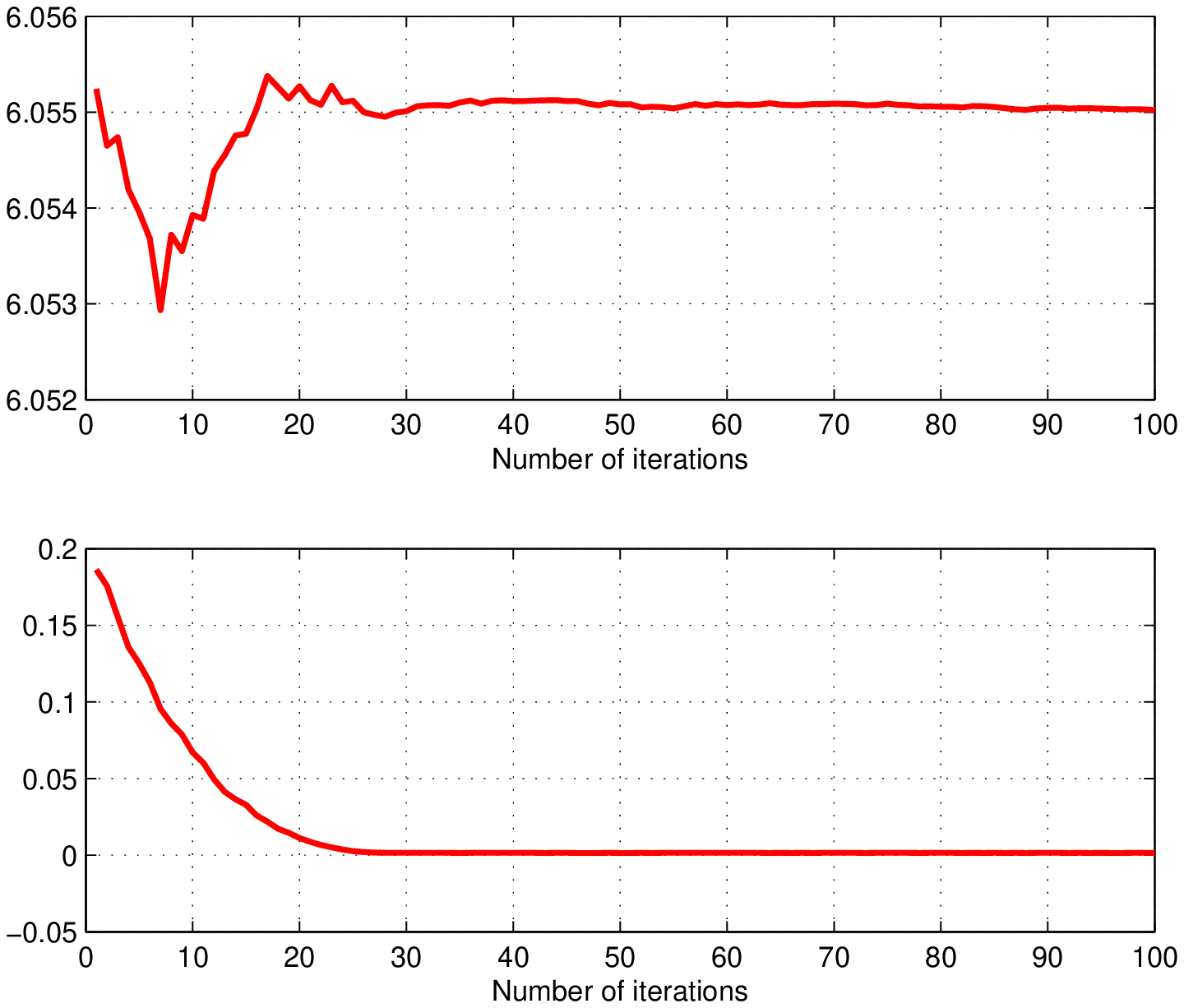}
\caption{Asian Option when $(\upsilon,K,\nu) = (0.1,45,\nu_\star)$. [left panel] successive
  directions of stratification $t \mapsto \mu^{(t)}$. $\mu^{(0)}$ is
  proportional to the vector $(1, \cdots, 1)$ so that the $d$ curves start from
  the same point $1/\sqrt{d}$. By convention, the first component of
  $\mu^{(t)}$ is positive. [top right] successive estimations of the quantity
  of interest $t \mapsto \mathcal{E}^{(t)}$. [bottom right] successive values
  of the variance $t \mapsto (\sum_{\bi} p_\bi \ \hat{\sigma}_\bi^{(t)})^2$;
  the limiting value is $0.002$.}
\label{fig:Asie:CvgMu}
\end{figure}

We run {\tt AdaptStr} when $(\upsilon,K) = (0.1,45)$: on
Figure~\ref{fig:Asie:directions}, the optimal drift vector $\nu_\star$, the
direction $\mu^{(N)}$ obtained after $N$ iterations of {\tt AdaptStr}, and the
directions of stratification $\mu_\g,\mu_{\mathrm{reg}},\mu_l$ are plotted.
 In Figure~\ref{fig:Asie:CvgMu}, the successive directions $t \mapsto
 \mu^{(t)}$, the successive estimations of the quantity of interest $t \mapsto
 \mathcal{E}^{(t)}$ and of the variance $t \mapsto (\sum_{\bi} p_\bi
 \hat{\sigma}_\bi^{(t)})^2$ are displayed.  We observe that $\{\mu^{(t)}, t
 \geq 0\}$ converges to the direction $\mu_\g$, and the convergence takes place
 after about 30 iterations.  We find the same pattern for a wide range of
 parameter values. The choice of the stratification direction has a major
 impact on the variance of the estimate $\mathcal{E}^{(t)}$ as shown on
 Figure~\ref{fig:Asie:CvgMu}~[bottom right].  Along the iterations of the
 algorithm, the variance decreases from $0.1862$ to $0.0016$. We also observed
 that the convergence of the algorithm and the limiting values were independent
 of the initial values $(\mu^{(0)}, M^{(0)})$ (these results are not reported
 for brevity).  These initial values (and the choice of the sequence
 $\{\gamma^{(t)}, t \geq 1 \}$) only influence the number of iterations
 required to converge.
{\tt AdaptStr} can also be seen as a procedure that computes a stratification
direction and provides the associated optimal allocation. These quantities can
then be used for running a (usual) stratification procedure with $M$ draws and
for the optimal allocation. By doing such with $M=20 \, 000$, we obtain an
estimate of the quantity $\PE[\phi(Y)]$ equal to $6.05$ and of the variance
equal to $0.002/M$. We can compare these results to the output of ${\tt GHS}$:
this yields the same estimator of $\PE[\phi(Y)]$ and a larger standard
deviation equal to $0.014/M$.  Observe that since $\mu^{(N)} = \mu_\g$, the
two algorithms differ from the allocations in the strata.

\begin{figure}[h]
  \centering \incg{200pt}{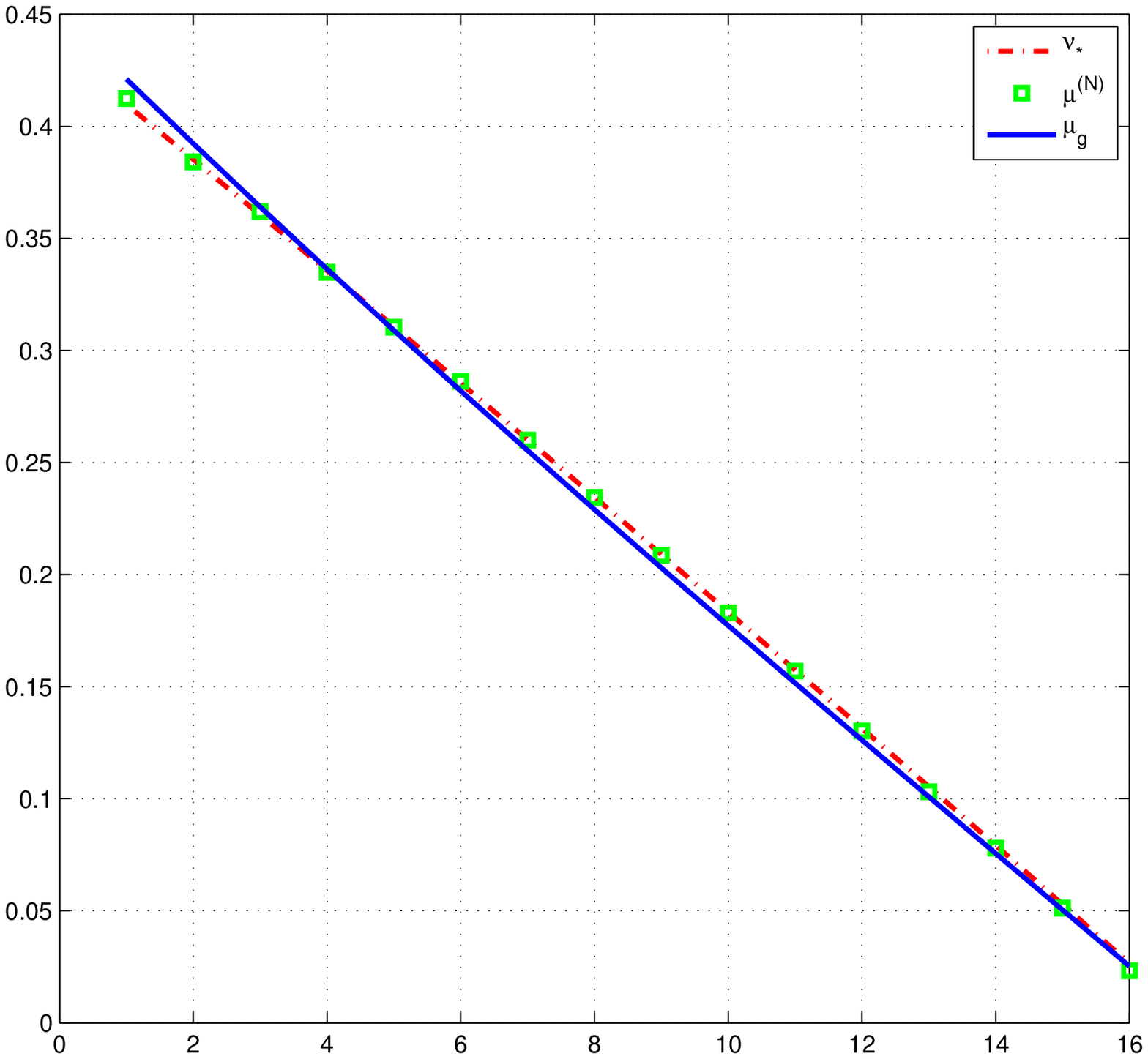}
  \incg{200pt}{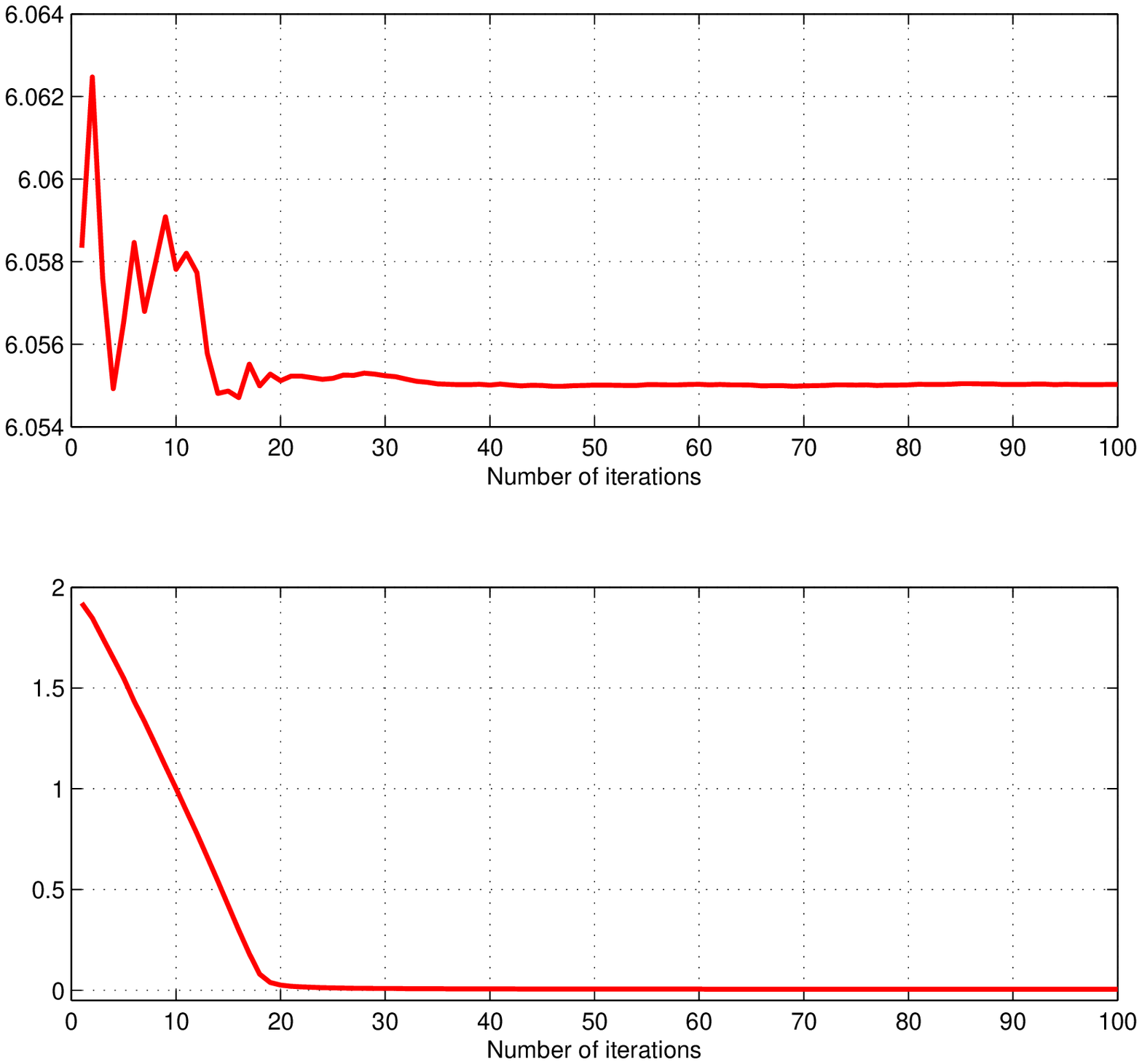}
\caption{Asian Option with $(\upsilon, K, \nu) = (0.1,45,0)$: [left panel] the limiting
  direction $\mu^{(N)}$ and for comparison, $\mu_g$ and $\nu_\star$ normalised
  to have norm $1$ ($\nu_\star \to \nu_\star/0.42$). [top right] successive
  estimations of the quantity of interest $t \mapsto \mathcal{E}^{(t)}$.
  [bottom right] successive values of the variance $t \mapsto (\sum_{\bi} p_\bi
  \ \hat{\sigma}_\bi^{(t)})^2$: the limiting value is $0.004$.}
\label{fig:Asie:DriftZero}
\end{figure}

We conclude this study of {\tt AdaptStr} by illustrating the role of the drift
vector $\nu$ (see Eq.~\ref{eq:Cameron-Martin}). We report on
Figure~\ref{fig:Asie:DriftZero} the limiting direction $\mu^{(N)}$, the
estimates $t \mapsto \mathcal{E}^{(t)}$ and the variance $t \mapsto (\sum_{\bi}
p_\bi \hat{\sigma}_\bi^{(t)})^2$ when $\nu=0$.
The limiting direction $\mu^{(N)}$ slightly differs from $\mu_\g$ and is close
to $\nu_\star$. Moreover, the variance reduction is weaker: the limiting value
of $t \mapsto (\sum_{\bi} p_\bi \hat{\sigma}_\bi^{(t)})^2$ is $0.004$.  The
efficiency of the adaptive stratification procedure {\tt AdaptStr} is thus
related to the drift vector $\nu$ in (\ref{eq:Cameron-Martin}); similar
conclusions are reached in \cite{glasserman:heidelberger:shahabuddin:1999} (see
also \cite{glasserman:2004}).

We report in Tables~\ref{tab:Asie} and \ref{tab:AsieLHS} the variance of
different estimators, as described in Section~\ref{sec:assessing_efficiency}.
\begin{center} {\em Insert
    Tables~\ref{tab:Asie} and \ref{tab:AsieLHS} about here}
\end{center}
Consider first the case $\nu=0$.  When the volatility of the asset is low
$\upsilon=0.1$ and the strike is in-the-money, the performance of the adaptive
stratification estimator ''{\tt AdaptStr}'' and of the stratification estimator
with fixed direction $\mu_{\mathrm{reg}}$, $\mu_\star$ and $\mu_l$ and with
optimal allocations are equivalent. We observe indeed that the directions
$\mu^{(N)}$,$\mu_{\mathrm{reg}}$, $\mu_\star$ and $\mu_l$ are almost colinear.
Compared to the plain Monte Carlo, the variance reduction factor is equal
to 2500.  The LH estimator with a rotation along any of the directions
$\mu^{(N)}$,$\mu_\star$ and $\mu_{\mathrm{reg}}$ outperforms all the stratified
estimators: the variance reduction is by a factor $10500$.  This reduction in
the variance strongly depends upon the choice of the rotation: the LH estimator
with no rotation implies a variance reduction by a factor $150$.

When the volatility of the asset is high $\upsilon=1$, the conclusions are
markedly different. Consider e.g. the case when the option is out of the
money ($K=65$).
The adaptive stratification estimator ''{\tt AdaptStr}'' provides a reduction
of variance by a factor $150$, which is again similar to the variance reduction
afforded by the stratification with fixed directions $\mu_{\mathrm{reg}}$,
$\mu_\star$ and $\mu_l$, and optimal allocation; {\tt AdaptStr} outperforms
stratified estimators with any of the fixed direction $\mu_{\mathrm{reg}}$,
$\mu_\star$ or $\mu_l$ by a factor $13$ when allocation is proportional.  The
LH estimator with no rotation only provides a reduction in variance by a factor
$1.7$; when the rotation along $\mu^{(N)}$ is applied, the reduction is by a
factor $65$. Here again, the LH estimator is very sensitive to the choice of
the orthogonal matrix $O$.

The use of the drift $\nu=\nu^\star$ improves the variance of all the
stratified estimators by a factor $2$ to $10$, depending on the choice of the
stratification direction; and by a factor $10$ to $25$ for ``{\tt MC}''. Here
again, ''{\tt AdaptStr}'' is the best stratified estimator; its performance can
be approached by stratification estimators with fixed directions, but the
choice of this fixed direction depends crucially upon the values of the
volatility and the strike. The vectors $\mu^{(N)}$ and $\mu_{\mathrm{reg}}$ are almost
colinear in many cases e.g. when $(\upsilon, K,\nu) = (0.1; 45; \nu_\star)$, but not always as observed in the case
$(\upsilon, K,\nu) = (1; 65;
\nu_\star)$ from the variances given in Table~\ref{tab:Asie}. It is interesting to note that the use of the drift $\nu_\star$
does not always improve the variance of the LH estimator.

These experiments show that the choice of the stratification direction and of
the allocation is crucial.  For example, in the case $(\upsilon,K,\nu)=
(0.5,65,0)$, the adaptive stratification estimator improves upon the
stratification estimator with fixed direction $(1,\cdots, 1) /\sqrt{d}$ and
optimal allocation by a factor $60$ (and by a factor $190$ when proportional
allocation is used) - these results are not reported in the tables for brevity
since this direction is rarely optimal - . Even if simple guesses for the
direction reduce the variance, this reduction can be improved (by a factor
$20$) when optimal allocation is used; this allocation is unknown and has to be
learnt.  In these examples, LH outperforms in many cases stratification
procedures provided it is applied with a rotation $O$: the rotation along
$\mu^{(N)}$ outperforms LH with no rotation and when compared with other simple
guess rotations, it provides similar or better variance reduction.  All these
remarks strongly support the use of adaptive procedures.

\subsection{Options with knock-out at expiration}
\begin{figure}[h]
\centering
\incg{200pt}{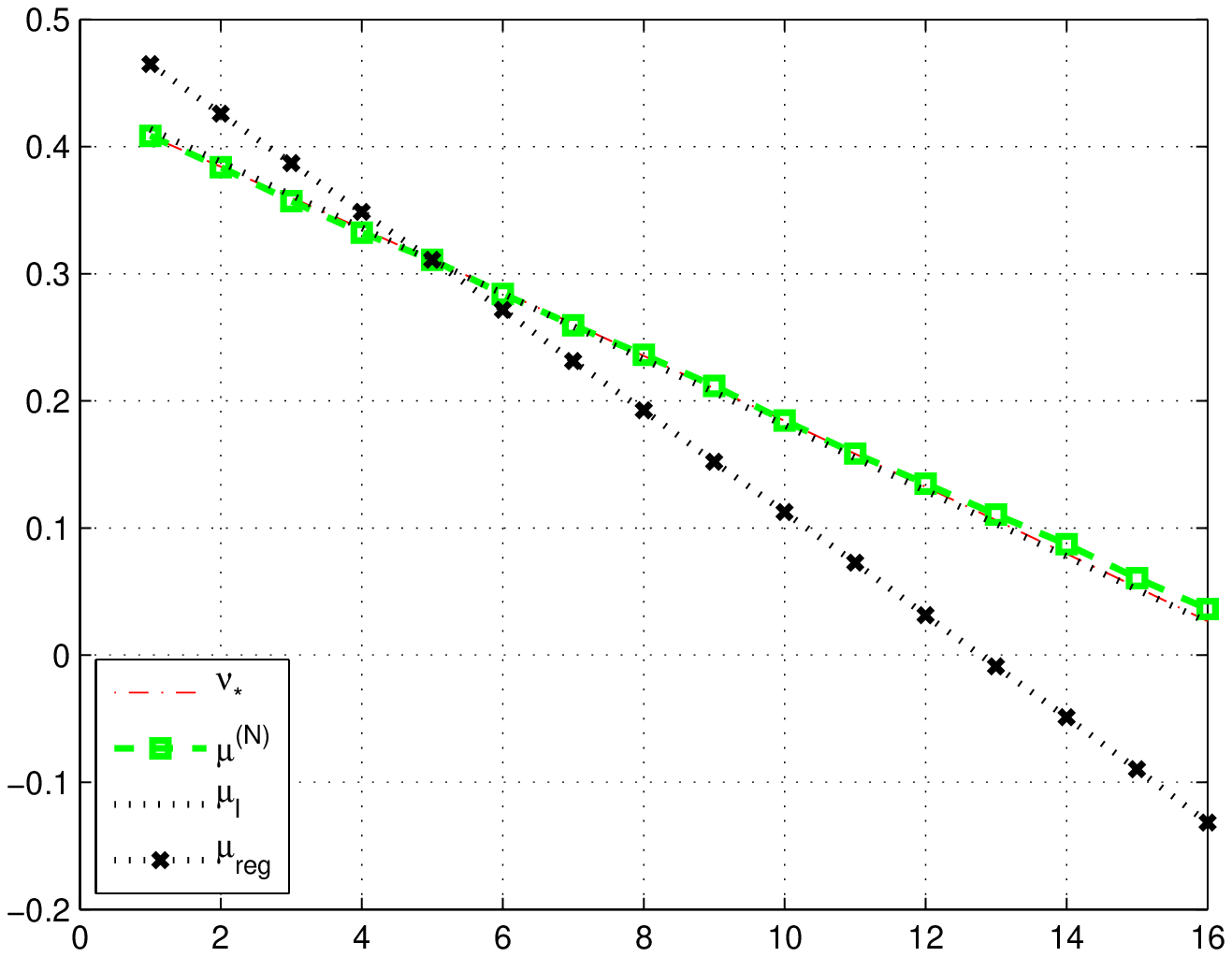} \incg{200pt}{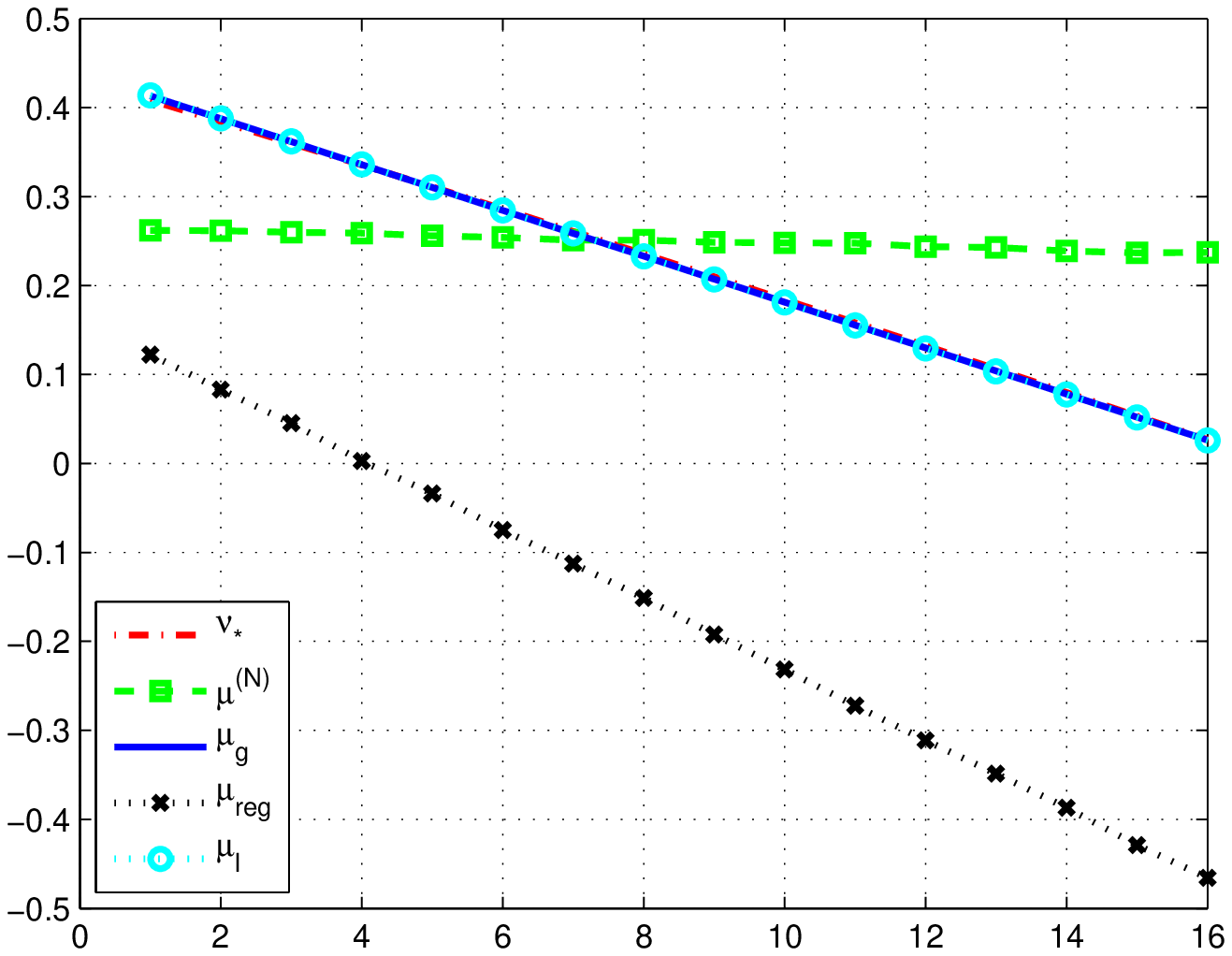}
\caption{Barrier Option when $(K,B) = (50,60)$ and [left panel]  $\nu =0$ or [right panel] $\nu = \nu_\star$: directions $\mu^{(N)}$,  $\mu_{\mathrm{reg}}$, $\mu_l$, $\nu_\star$ and $\mu_g$.  $\nu_\star$ has been scaled to have norm $1$ ($\nu_\star \leftarrow \nu_\star /
  0.84$)}
 \label{fig:Barrier:directions}
\end{figure}

A knock-out barrier option is a path-dependent option that expires worthless if
the underlying reaches a specified barrier level.  The payoff of this option is
given by
\[
\Xi(y) = \exp(-rT) \left( \frac{s_0}{d} \sum_{k=1}^d \exp\left((r - 0.5
    \sigma^2) \frac{k T}{d} + \sigma \sqrt{\frac{T}{d}} \sum_{j=1}^k y_j
  \right) -K \right)_+ \1_{\{S_T(y) \leq B\}} \eqsp,
\]
where $K$ is the strike price, $B$ is the barrier and $S_T(y)$ is the underlier
price modeled as
\[
S_T(y) = s_0 \exp\left((r - 0.5 \sigma^2) T + \sigma \sqrt{\frac{T}{d}}
  \sum_{j=1}^d y_j \right) \eqsp.
\]
In the numerical applications, we set $s_0 = 50$, $r=0.05$, $T=1$, $\sigma =
0.1$, $d=16$ and $(K,B) \in \{( 50,60), ( 50,80)\}$. We choose $\mu_l \propto
(d,d-1,\cdots, 1)$.  On Figure~\ref{fig:Barrier:directions}[left panel], we plot
$\mu^{(N)}$ in the case $(K,B,\nu) = (50,60,0)$; the limiting direction is
$\mu_l$.  On Figure~\ref{fig:Barrier:directions}[right panel], we plot $\mu^{(N)}$ in
the case $(K,B,\nu) = (50,60,\nu_\star)$; for comparison, we also plot
$\nu_\star$, $\mu_g$, $\mu_{\mathrm{reg}}$ and $\mu_l$. This is an example
where the optimal stratification direction $\mu^{(N)}$ does not coincide with
the different directions of stratification ($\mu_\g$, $\mu_{\mathrm{reg}}$ and
$\mu_l$); in this example, $\mu_g \sim \mu_l$ but the optimal direction of
stratification is close to $(1, \cdots, 1)/\sqrt{d}$.  The direction associated
to the regression estimator is far from being optimal.

We report in Tables~\ref{tab:Barrier} and \ref{tab:BarrierLHS} the variances of
the different estimators, as described in
Section~\ref{sec:assessing_efficiency}.
\begin{center}
  {\em Insert Table~\ref{tab:Barrier} and \ref{tab:BarrierLHS} about here}
\end{center}
Consider first the case where the drift $\nu$ is set to $0$. When $(K,B)=(50,60)$
(the option is at the money, and the barrier is close to the money), the
adaptive stratification provides a variance reduction by a factor 10 with
respect to the plain Monte Carlo estimator.  In this case, the stratification
directions $\mu_l$ and $\mu_\star$ (with an optimal allocation) perform almost
as well (and $\mu_l$ and $\mu_\star$ are close to $\mu^{(N)}$ at the
convergence), while $\mu_{\mathrm{reg}}$ provides a higher variance.  For the LH estimator, the variance reduction is only
by a factor $1.5$. It is worthwhile to note that the best choices for the rotation
, $\mu^{(N)}$ and $\mu_\star$, lead to a variance thrice the one of ''{\tt AdapStr}''.  The use of the drift vector $\nu_\star$ improves the variance
of the adaptive stratification by a factor $1.8$; the optimal stratification
vector is now $\mu^{(N)}$ which surpasses $\mu_{\mathrm{reg}}$.  The variance
of the LH estimator is also reduced.

When $(K,B)=(50,80)$ (the barrier is out of the money), a factor reduction 2800
is obtained by the adaptive stratification estimator ''{\tt AdaptStr}''; a
similar variance reduction is achieved using the stratified estimator with
direction $\mu_{\mathrm{reg}}$ and with optimal allocation. For the LH
estimator, the variance reduction is by a factor 2000, when the rotation is
$\mu_{\mathrm{reg}}$.  The use of the drift vector $\nu_\star$ improves the
behavior of all the algorithms: the variance of ``{\tt AdaptStr}'' is reduced
by a factor $3.8$. Finally, LH with rotation $\mu^{(N)}$ reduces the variance
of LH with no rotation by a factor $1200$.

To conclude, this example shows again the interest of adaptive procedures in
order to find a stratification direction, the optimal allocation or a rotation
in LH.

\subsection{Basket options}
Consider a portfolio consisting of $d$ assets.  The portfolio contains a
proportion $\alpha_k$ of asset $k$, $k \in \{1, \dots, d \}$. The price of each
asset is described by a geometric Brownian motion (under the risk neutral
probability measure)
\[
\frac{dS_t^{(k)}}{S_t^{(k)}} = r \, dt + \upsilon_{k} \, d W_t^{(k)}
\]
but the standard Brownian motions $\{W^{(k)}_{.}, k \in \{1, \dots, d\} \}$ are
not necessarily independent. For any $t\geq s$ and $k \in \{1, \dots, d \}$
\[
\ln S^{(k)}_t = \ln S^{(k)}_s + \left(r - 0.5 \upsilon^2_k\right) (t-s) +
\upsilon_k \sqrt{t-s} \tilde Y_k
\]
where $\tilde Y = (\tilde Y_1, \dots, \tilde Y_d) \sim \mathcal{N}_d(0,
\Sigma)$. The $d \times d$ matrix $\Sigma$ is a positive semidefinite matrix
with diagonal coefficients equal to $1$. Therefore, the variance of the
log-return on asset $k$ in the time interval $[s,t]$ is $(t-s) \upsilon^2_{k}$,
and the covariance between the log-returns $i,j$ is $(t-s) \upsilon_i
\upsilon_j \Sigma_{i,j}$. It follows that $\Sigma_{i,j}$ is the correlation
between the log-returns.  The price at time $0$ of a European call option with
strike price $K$ and exercise time $T$ is given by $\PE[\Xi(Y)]$ where
\[
\Xi(y) = \exp(-rT) \left( \sum_{k=1}^d \alpha_k s_0^{(k)} \exp\left((r - 0.5
    \upsilon^2_k) T + \upsilon_k \sqrt{T} \tilde y_k\right) -K \right)_+
\]
and $\tilde y = \sqrt{\Sigma} y$ ($\sqrt{\Sigma}$ denotes a square root of the
matrix $\Sigma$ i.e. solves $MM^T = \Sigma$). In the numerical applications,
$\Sigma$ is chosen to be $\Sigma_{i,j} =1_{\{i=j\}}+c1_{\{i\neq j\}}$, $\alpha_k =
1/d$, $r=0.05$, $T=1$, and $d=40$.  We consider $(c,K) \in \{(0.1, 45), (0.5,
45),(0.9, 45) \}$. The initial values $\{s_0^k, k \leq d\}$ are drawn from the 
uniform distribution in the range $[20, 80]$; the volatilities $\{\upsilon_k, k
\leq d\}$ are chosen linearly equally spaced in the set $[0.1, 0.4]$. The
assets are sorted so that $\upsilon_1 \leq \cdots \leq \upsilon_d$. We choose
\[
\mu_l \propto \left(\alpha_1 s_0^{(1)} \exp((r - 0.5 \upsilon^2_1) T) \upsilon_1, \cdots, \alpha_d s_0^{(d)} \exp((r - 0.5 \upsilon^2_d) T) \upsilon_d \right) \sqrt{\Sigma} \eqsp.
\]
In the case $(c,K,\nu) = (0.1,45,\nu_\star)$, we plot on
Figure~\ref{fig:Asie:directions}[right panel] the limiting direction $\mu^{(N)}$ and
for comparison, the directions $\nu_\star$, $\mu_g$, $\mu_{\mathrm{reg}}$,
$\mu_l$. We report in Tables~\ref{tab:Basket} and \ref{tab:BasketLHS} the variance of
the different estimators, as described in
Section~\ref{sec:assessing_efficiency}.
\begin{center}
  {\em Insert Tables~\ref{tab:Basket} and \ref{tab:BasketLHS} about here}
\end{center}
In this example again the adaptive stratification estimator improves upon the
best stratified estimator with (non-adaptive) stratification direction and optimal
allocation.  Here again, the optimal allocation improves the variance reduction,
by a factor $15$ for example in the case $(c,K,\nu) = (0.9, 45, 0)$ for the
fixed directions $\mu_{\mathrm{reg}}$, $\mu_\star$ or $\mu_l$.  It is
interesting to note that the variance reduction with respect to the plain Monte
Carlo using "{\tt AdaptStr}'' ranges from 100
($c=0.1$,$K=45$) to 2500 ($c=0.9$,$K=45$) whereas the use of the drift
$\nu^\star$ allows only a reduction by a factor 10. The choice of the
stratification direction plays a more important role than the choice of drift direction.

The comparison with the LH estimator is more difficult, because this estimator
behaves totally differently from the stratified estimator. First, except
for $c=0.9$, the use of
the drift $\nu_\star$ increases the variance: whereas the effect of the drift is always
markedly beneficial for the Monte Carlo estimators, the drift $\nu_\star$ may increase the
variance by a factor as large as 15 (when $(c,K)=(0.1,45)$).
Second, LH with rotation always outperforms (adaptive) stratification: the
main difficulty stems from the choice of the rotation but the obtained results show that
rotation along $\mu^{(N)}$ provides either the maximal variance reduction or 
variance reduction similar to the best rotation among the three considered.
Finally, the performance of the estimator is extremely sensitive to the choice
of the simulation setting: when the correlation among the assets is large
$c=0.9$, the choice of the first vector of the orthogonal matrix $O$ becomes
crucial. With drift $\nu_\star$, rotation along $\mu^{(N)}$ (resp. $\mu_{\mathrm{reg}}$) improves LH
with no rotation by a factor $22500$ (resp. $3000$).

\subsection{Heston model of stochastic volatility}
We consider a last example which  is not covered by the methodology
presented in \cite{glasserman:heidelberger:shahabuddin:1999}. We price an Asian option in the Heston model, 
specified as follows
$$d\xi_t=k(\theta-\xi_t)dt+\sigma\sqrt{\xi_t}\,dW^1_t,\;\;dS_t=
rS_tdt+\sqrt{\xi_t}\,S_t\,(\rho dW^1_t+ \sqrt{1-\rho^2}dW^2_t),\;\;dX_t= S_tdt$$
where $\{W^1_t, t\geq 0\}$ and $\{W^2_t, t\geq 0\}$ are two independent
Brownian motions, $r$ is the risk free rate, $\sigma>0$
the volatility of the volatility process $\xi$, $k\geq 0$ the mean reversion rate, $\theta\geq 0$ the long run average
volatility, and $\rho\in[-1,1]$ the correlation rate. The price
of an Asian Call option with strike $K$ and maturity $T$ is
\begin{equation}
\label{pricehes}
\PE \left[\exp(-rT)\left(\frac{1}{T}X_T-K \right)_+\right].
\end{equation}
In our tests, we have chosen the parameters so that $\sigma^2\leq 4 k\theta$.
This enabled us to replace by Gaussian increments the finitely-valued random
variables used to discretize $W^1$ in the scheme proposed in
\cite{alfonsi:2008} to approximate the SDE satisfied by $(S_t,X_t,\xi_t)$. We
refer to \cite{alfonsi:2008} for a precise description of the scheme that we
used.  The resulting approximation $\hat{X}_d$ of $X_T$ is generated from a
vector $Y=(Y_1,\ldots,Y_d,Y_{d+1},\ldots,Y_{2d})\sim\mathcal{N}_{2d}(0,\Id)$
corresponding to the increments of $(W^1,W^2)$ and a vector
$B=(B_1,\ldots,B_d)$ of independent Bernoulli random variables with parameter
$0.5$. The price \eqref{pricehes} is then approximated by
$$\PE [\exp(-rT)\left(
  \hat{X}_d-K \right)_+].$$



In the following tests we keep $\nu$ in (\ref{eq:Cameron-Martin}) equal
to zero and do not stratify the random vector $B$.  

We choose  $S_0=100$, $\theta=0.01$, $k=2$, $\sigma=0.2$, $T=1$, $r = 0.095$,
$\rho=-0.5$ and $(\xi_0, K) \in \{(0.01, 120), (0.01, 100), (0.01,80),
(0.04,130), (0.04,100), (0.04,70) \}$. The discretization step of the scheme is
$d=50$.

\begin{figure}[h]
  \centering \incg{200pt}{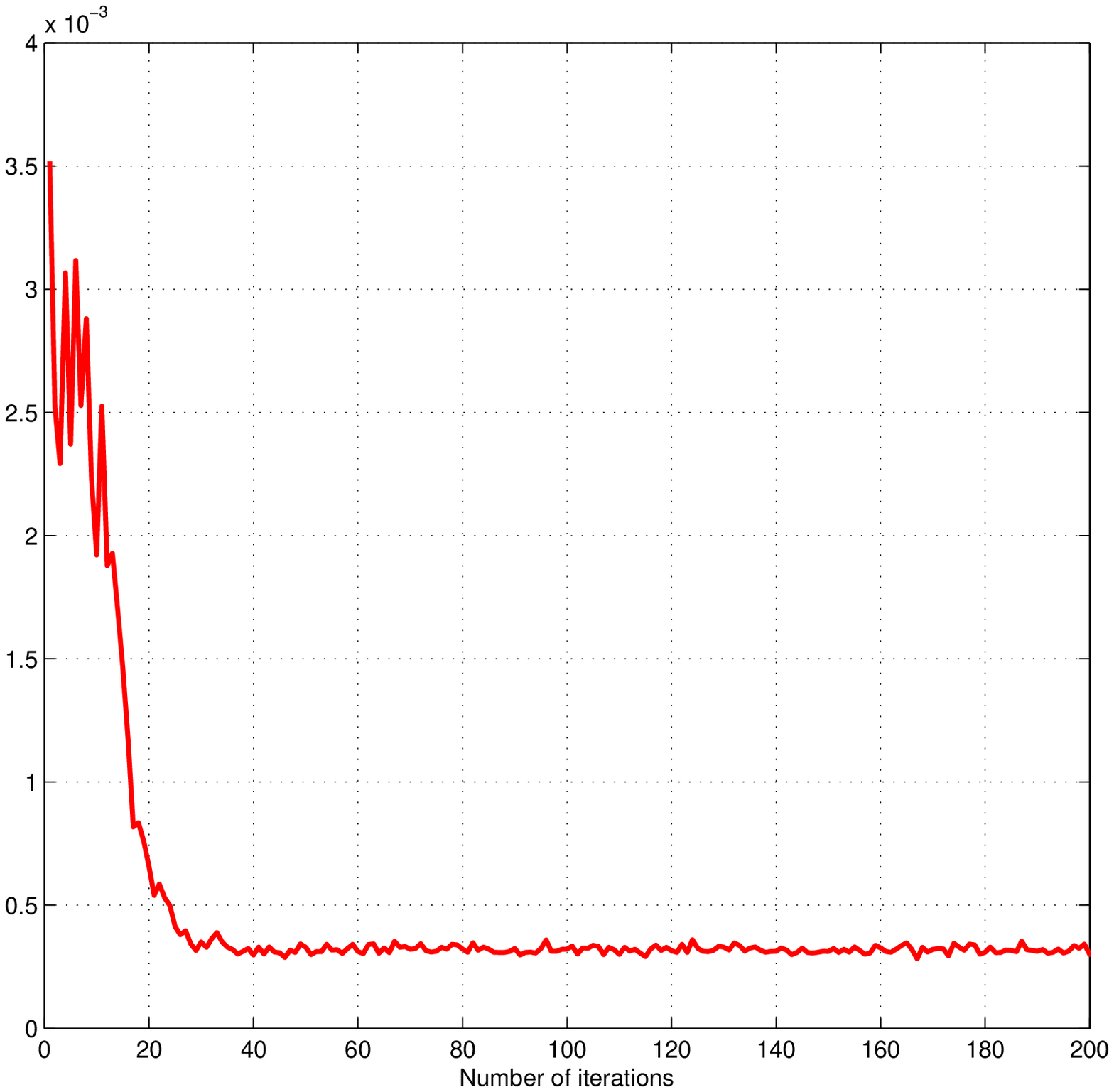} \incg{200pt}{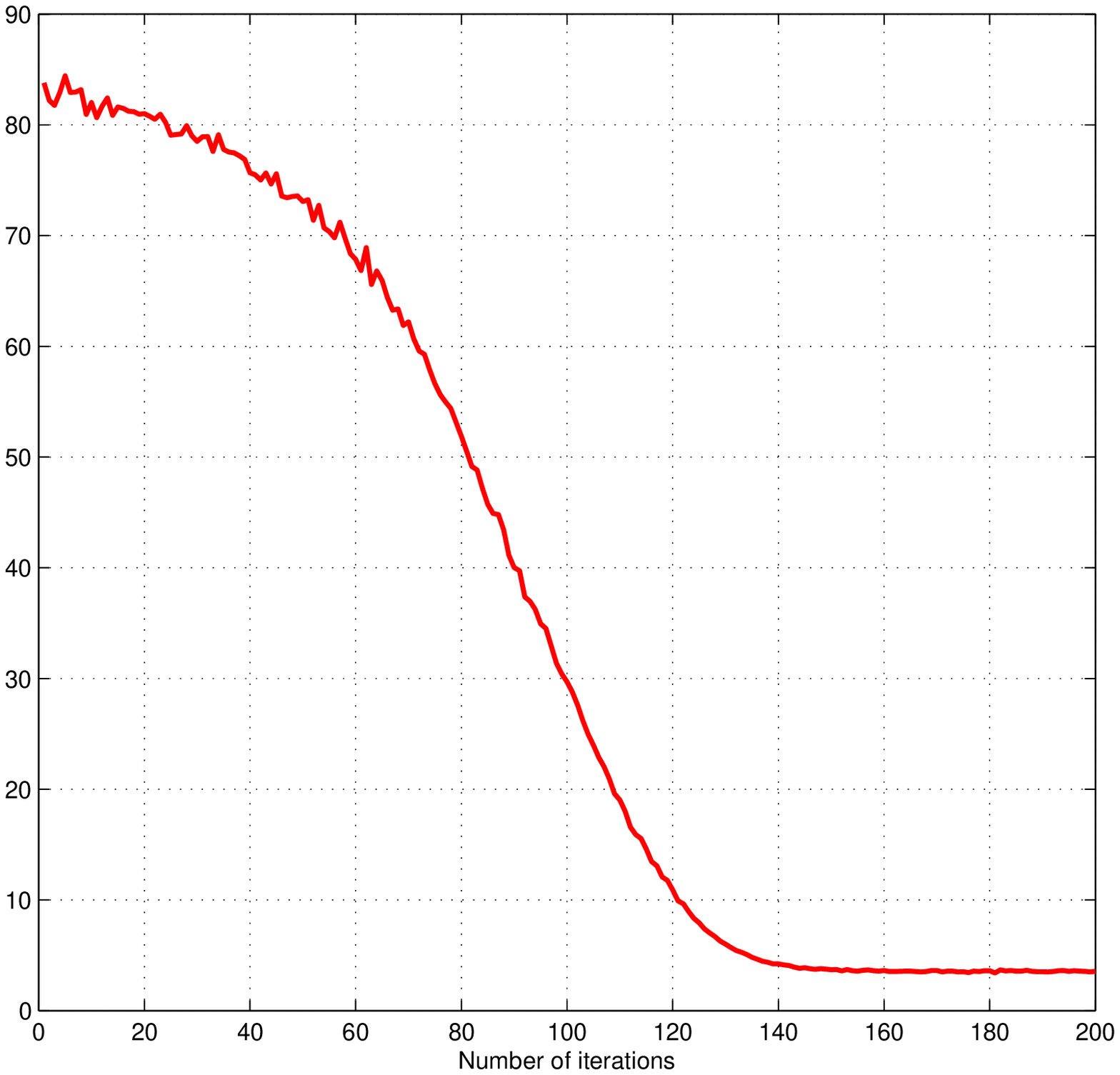}

\caption{Asian Call Option in Heston model: successive estimations of 
  the variance $t \mapsto (\sum_{\bi} p_\bi \ \hat{\sigma}_\bi^{(t)})^2$ [left
  panel] when $(\xi_0, K) = (0.01; 120)$; [right panel] when $(\xi_0, K) =
  (0.04; 70)$}
\label{fig:evolcost:heston}
\end{figure}

On Figure~\ref{fig:evolcost:heston} we plot the successive estimations of the
variance $t \mapsto (\sum_{\bi} p_\bi \hat{\sigma}_\bi^{(t)})^2$, when $(\xi_0,
K) = (0.01; 120)$ and $(\xi_0, K) = (0.04; 70)$.

\begin{figure}[h]
   \begin{center}
   \incg{200pt}{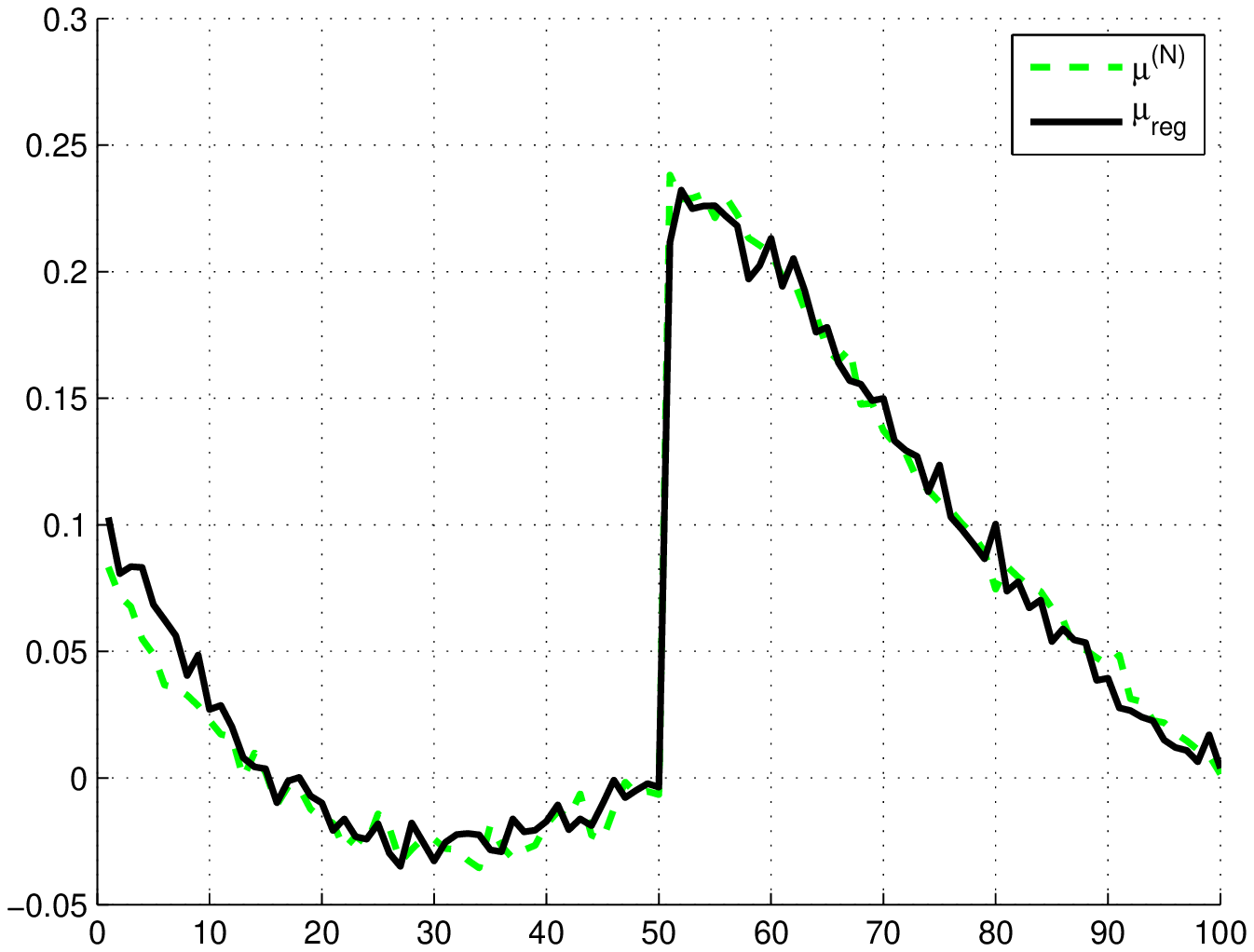}   \incg{200pt}{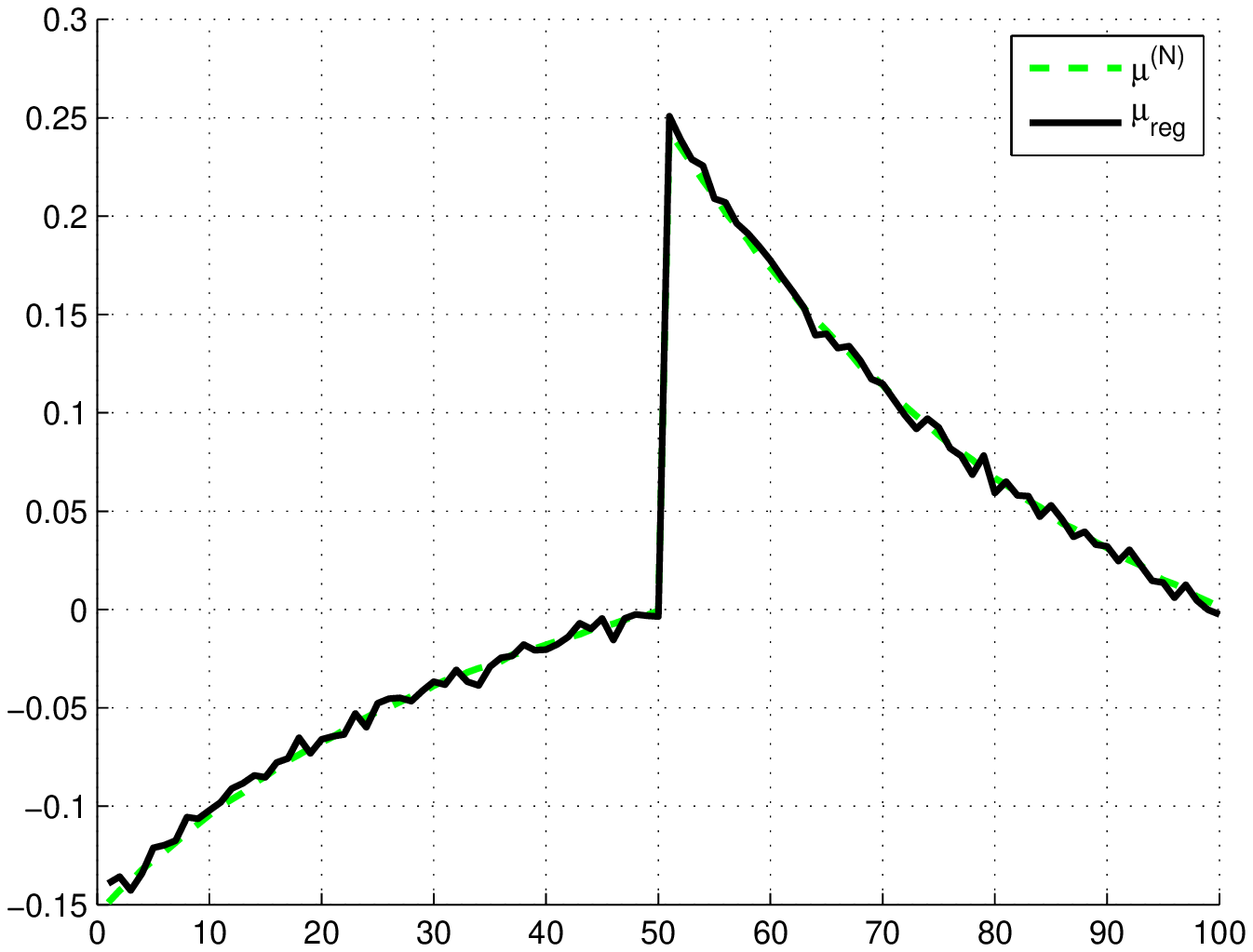}
\end{center}
\caption{Asian Call Option in Heston model: vectors $\mu^{(N)}$ and  $\mu_{\mathrm{reg}}$ in the cases [left panel] $(\xi_0, K) = (0.01; 120)$; [right panel]  $(\xi_0, K) = (0.04; 70)$}
\label{fig:dirfin:heston}
\end{figure}

We plot on Figure \ref{fig:dirfin:heston} the components of $\mu^{(N)}$ with
respect to the component index in the cases $(\xi_0, K) = (0.01; 120)$ and
$(\xi_0, K) = (0.04; 70)$.

We report in Tables~\ref{tab:Heston} and \ref{tab:HestonLHS} the variances of
some estimators described in Section~\ref{sec:assessing_efficiency}.
\begin{center}
  {\em Insert Table~\ref{tab:Heston} and \ref{tab:HestonLHS} about here}
\end{center}

The first observation is that even in this case, {\tt AdaptStr} still works and
provides variance reduction when compared to Monte Carlo. It is all the more
efficient than the option is out of the money : when $(\xi_0, K,\nu) = (0.01,
120, 0)$, the variance reduction is by a factor $85$; when $(\xi_0, K,\nu) =
(0.04, 130, 0)$, the variance reduction is by a factor $105$. {\tt AdaptStr}
and stratification with fixed direction $\mu_{\mathrm{reg}}$ are equivalent,
provided the last one is applied with optimal allocation, just necessitating
again iterative procedures.

We can wonder on the effect of the moneyness and the volatility of the model on
the variance reduction. As shown in Table \ref{tab:Heston}, in general
the achieved variance reduction is larger when the option is out of the money
(for $K=130$ and $\xi_0=0.04$ the variance is divided by nearly $105$ when
using {\tt AdaptStr}). We also observe that stratification procedures
outperform LH samplers when the option is out of the money, but LH is
equivalent to stratification when the option is in the money.

\section{Proofs}
\label{sec:Proofs}
\subsection{Proofs of Section~\ref{sec:AsymptoticRegime}}
\label{sec:Proof:TheoryMultiDim}
In the sequel, we denote $\mathcal{I}_m \eqdef \{1, \dots, I\}^m$ and $p_\bi \sigma_\bi \eqdef \sqrt{\left(\int_{\strata{\bi}}\marg \ d\leb \right)
  \left(\int_{\strata{\bi}}\mmarg\marg \ d\leb
  \right)-\left(\int_{\strata{\bi}}\smarg\marg \ d\leb \right)^2}$ in
place of $p_\bi(\mu) \sigma_\bi(\mu)$.
 \begin{lemma}\label{arrondis} Let $\{\strata{\bi}, \bi \in \mathcal{I}_m
   \}$ be given by (\ref{eq:definition-stratum}).
\begin{enumerate}[(i)]
\item \label{arrondis1}
$\forall \epsilon >0,\;\forall M > \epsilon^{-1},\;\sup_{\Q: \inf_{\bi \in \mathcal{I}_m} q_\bi \geq \epsilon } \ \left|M
  \varsigma_{I, M}^2(\mu,g,\Q)-\sum_{\bi \in \mathcal{I}_m}
  \frac{p_{\bi}^2\sigma_{\bi}^2}{q_{\bi}} \right| \leq \frac{1}{M \epsilon (\epsilon-M^{-1})}
\ \PVar[\phi(Y)].$
\item \label{arrondis2} Assume that $\essinf{g \cdot \leb}{\chi /g}>0$ and
    $\esssup{\chi \cdot \leb}{\marg/\chi}<+\infty$.  Let $\epsilon>0$.
  For any $(I,M)$ such that $M I^{-m} \essinf{g \cdot \leb}{\chi /g} \geq
  1+\epsilon$
$$\left|M \varsigma_{I, M}^2(\mu,g,\Q_\chi)-\sum_{\bi \in \mathcal{I}_m}\frac{p_{\bi}^2\sigma_{\bi}^2}{q_{\bi}(\chi)} \right| \leq
\frac{(1 + \epsilon^{-1}) \PVar[\phi(Y)]}{ \essinf{g \cdot \leb}{\chi /g}}
\ \frac{I^m}{M} \ \left(\esssup{\chi \cdot \leb}{\frac{\marg}{\chi}} \wedge
\frac{I^m}{\essinf{g \cdot \leb}{\chi /g}}\right).$$
\item \label{arrondis3} For any positive integers $M,I$ and real $\epsilon>1$,
\[
\left|M \varsigma_{I, M}^2(\mu,g,\Qopt(\mu))-\sum_{\bi \in \mathcal{I}_m} \frac{p_{\bi}^2
  \,\sigma_{\bi}^2 }{\qopt{\bi}(\mu)}\right| \leq \PVar[\phi(Y)] \ \left(
  (1+\epsilon)\frac{I^m}{M}+ \frac{1}{\epsilon-1}\right) \eqsp,
\]
where $\Qopt(\mu) = \{ \qopt{\bi}(\mu), \bi \in \mathcal{I}_m \}$ is the optimal allocation defined by
(\ref{eq:optimal-allocation-qi}).
   \end{enumerate}
 \end{lemma}
\begin{proof}
By definition of $M_\bi$ (see Eq. \ref{eq:definition-number-draws}),
  $M_\bi=0$ when $q_\bi=0$ and $M_\bi
  \geq 1$ when $q_\bi \geq M^{-1}$. One may have $M_\bi=1$ when
  $q_\bi\in(0,M^{-1})$ but then $MM_\bi^{-1}\leq q_\bi^{-1}$. Hence,
 \begin{equation}  \label{proof:arrondis}
   \left| M \; \varsigma^2_{I,M}(\mu,g,\Q) - \sum_{\bi \in \mathcal{I}_m: q_\bi >0} \ \frac{p_{\bi}^2\sigma_{\bi}^2}{q_{\bi}} \right|
   \\
   \leq \sum_{\bi \in \mathcal{I}_m, q_\bi \geq 1/M}\ \left| \frac{M
       q_{\bi} - M_{\bi}}{ M_{\bi}} \right|
   \frac{p_{\bi}^2\sigma_{\bi}^2}{q_{\bi}} +\sum_{\bi \in \mathcal{I}_m:
     0<q_\bi < 1/M} \ \frac{p_{\bi}^2\sigma_{\bi}^2}{q_{\bi}}
 \end{equation}
      \textit{(i)} When $\inf_{i\in\mathcal{I}_m} q_\bi \geq \epsilon > M^{-1}$, the second term in the
      rhs is null and since by \eqref{eq:definition-number-draws},
      $Mq_\bi-1<M_\bi<Mq_\bi+1$,
\begin{equation}\label{majot1}
   \sum_{\bi \in \mathcal{I}_m, q_\bi \geq 1/M}\ \left| \frac{M
       q_{\bi} - M_{\bi}}{ M_{\bi}} \right|
   \frac{p_{\bi}^2\sigma_{\bi}^2}{q_{\bi}}\leq M^{-1} \ \left(\sup_{\bi \in \mathcal{I}_m, q_\bi \geq 1/M} p_\bi
  q_\bi^{-1} \right) \sum_{\bi \in \mathcal{I}_m, q_\bi \geq 1/M}\ (q_\bi -
M^{-1})^{-1} \ p_\bi \sigma_\bi^2,
\end{equation}
which yields the desired result upon noting that $ p_\bi q_\bi^{-1}\leq
q_\bi^{-1} \leq \epsilon^{-1}$ and $\sum_{\bi} p_\bi \sigma_\bi^2 \leq
\PVar[\phi(Y)]$.

\textit{(ii)} Under the stated assumptions, $ q_\bi(\chi) = \int_{\strata{\bi}}
\chi d \leb \geq \essinf{g \cdot \leb}{\chi /g} \ \ I^{-m}$. Hence $M q_\bi
\geq 1 + \epsilon$ which implies that the second term in the rhs of
(\ref{proof:arrondis}) is null. This also implies that
$q_\bi - M^{-1} \geq \left(1 - \frac{1}{1+\epsilon} \right) \essinf{g \cdot
  \leb}{\chi /g} \ I^{-m}$.
We conclude the proof by combining this bound with \eqref{majot1} and the
following one :
\[
\frac{p_{\bi}}{q_{\bi}(\chi)}=\frac{\int_{\strata{\bi}} \marg \
  d\leb}{\int_{\strata{\bi}}\chi \ d\leb}\leq \esssup{\chi \cdot
  \leb}{\marg/\chi} \wedge \frac{1}{q_\bi(\chi)} \leq \esssup{\chi
  \cdot \leb}{\marg/\chi} \wedge \frac{I^m}{\essinf{g \cdot \leb}{\chi
    /g}}\eqsp.
\]

\textit{(iii)} Note that by convention, $p_\bi^2 \sigma_\bi^2 /
\qopt{\bi}(\mu) =0$ when $\qopt{\bi}(\mu) =0$.  By
definition of the optimal allocation (see Eq.~\ref{eq:optimal-allocation-qi}),
\[
p_\bi^2 \sigma_\bi^2/\qopt{\bi}(\mu)=
\qopt{\bi}(\mu) \left( \sum_{\bj} p_\bj \sigma_\bj\right)^2 \leq
\qopt{\bi}(\mu) \ \PVar[\phi(Y)] \eqsp.
\]
The second term in the rhs of (\ref{proof:arrondis}) is upper bounded by $I^m
M^{-1} \ \PVar[\phi(Y)]$. For the first term,
\begin{multline*}
  [ \PVar[\phi(Y)]]^{-1} \  \sum_{\bi \in \mathcal{I}_n,
    \qopt{\bi}(\mu) \geq 1/M}\ \left| \frac{M \qopt{\bi}(\mu) - M_{\bi}}{
      M_{\bi}} \right| \frac{p_{\bi}^2\sigma_{\bi}^2}{
    \qopt{\bi}(\mu) } \\
 \leq \sum_{\bi \in \mathcal{I}_m, 1/M \leq \qopt{\bi}(\mu)
    \leq \epsilon/M}\ \left| \frac{M \qopt{\bi}(\mu) - M_{\bi}}{ M_{\bi}} \right|
\   \qopt{\bi}(\mu)  + \sum_{\bi \in \mathcal{I}_m, \qopt{\bi}(\mu) \geq
    \epsilon/M}\ \left| \frac{M \qopt{\bi}(\mu) - M_{\bi}}{ M_{\bi}} \right|  \
  \qopt{\bi}(\mu) \eqsp.
\end{multline*}
For all $\bi$ such that $\qopt{\bi}(\mu) \geq 1/M$, $M_\bi^{-1} |M
\qopt{\bi}(\mu) - M_{\bi}| \leq 1$ which implies that
\[
\sum_{\bi \in \mathcal{I}_m, 1/M \leq \qopt{\bi}(\mu) \leq
  \epsilon/M}\ \left| \frac{M \qopt{\bi}(\mu) - M_{\bi}}{ M_{\bi}} \right|
\qopt{\bi}(\mu) \leq \frac{\epsilon I^m}{M}\eqsp.
\]
For all $\bi$ such that $\qopt{\bi}(\mu) \geq \epsilon/M$,
$ M_\bi^{-1} |M \qopt{\bi}(\mu) - M_{\bi}| \leq M_\bi^{-1} \leq (M
\qopt{\bi}(\mu) -1)^{-1} \leq (\epsilon -1)^{-1}$
which implies that
\[
\sum_{\bi \in \mathcal{I}_m, \qopt{\bi}(\mu) \geq \epsilon/M}\
\left| \frac{M \qopt{\bi}(\mu) - M_{\bi}}{ M_{\bi}} \right|
\qopt{\bi}(\mu) \leq (\epsilon -1)^{-1} \eqsp.
\]
\end{proof}

\begin{proof}{\emph{of Proposition~\ref{prop:multVarLimGalCase}}}
To prove the Proposition~\ref{prop:multVarLimGalCase}, we need the two
following Lemmas. The
first is a standard change of variables formula (see for example, \cite[Theorem
4.1.11]{dudley:2004}). Define
$G^{-1}(x_1, \dots, x_m) \eqdef (G_1^{-1}(x_1), \dots ,G_m^{-1}(x_m))$ where $G_k$ is the c.d.f. associated to the density $g_k$ on $\rset$.
\begin{lemma}
Let $h:\rset^m\rightarrow\rset$ be
  a measurable function. Assume that $h$ is nonnegative or is such that
  $\int_{\rset^m}|h| 1_{\{g>0\}} \ d\leb <+\infty$. Then, for all $0\leq
  v_k\leq w_k\leq 1$, $k \in \{1,\dots, I\}$
\begin{equation}
\label{chgtvar}
\int_{\prod_{k=1}^m [G^{-1}_k(v_k), G^{-1}_k(w_k)]}h 1_{\{g>0\}} \ d \leb=\int_{\prod_{k=1}^m [v_k, w_k]}\frac{h}{g} \circ G^{-1} \ d\leb \eqsp.
\end{equation}
\end{lemma}
The second technical Lemma is our key approximation result.
\begin{lemma}
\label{lem:ControleRestes}
Let $h,\gamma : \rset^m \to \rset$ be functions such that $\int_{\rset^m} \left( h^2+\gamma^2\right) /g \ d\leb <+\infty$.
Define for $\bi \in \mathcal{I}_m$,
\begin{equation}
  \label{eq:DefiResteR}
  R_{\bi}[h,\gamma]\eqdef
\int_{\strata{\bi}} h\gamma /g \ d\leb-I^m
\left(\int_{\strata{\bi}}h \ d\leb \right)
\left(\int_{\strata{\bi}}\gamma \ d\leb \right)
\eqsp.
\end{equation}
Then $\lim_{I\rightarrow +\infty}\sum_{\bi \in \mathcal{I}_m} |R_{\bi}[h,\gamma]|=0$.
\end{lemma}
\begin{proof}
  By polarization, it is enough to prove the result when $\gamma=h$ with
  $\int_{\rset^m}h^2 /g \ d\leb<+\infty$. This integrability condition
  ensures that $\leb$-\pp, $g=0$ implies $h=0$ and by \eqref{chgtvar}, one has
  $$R_{\bi}[h,h]=\int_{\prod_{k=1}^m  [(i_k-1)/I, i_k/I]}\frac{h^2}{g^2}\circ G^{-1} \
  d\leb-I^m \ \left(\int_{ \prod_{k=1}^m  [(i_k-1)/I, i_k/I] }\frac{h}{g}\circ G^{-1} \
    d\leb\right)^2 \eqsp,$$
  where the right-hand-side is non-negative by
  Cauchy-Schwarz inequality.  Set  $\tilde h(u) \eqdef \frac{h}{g}(G^{-1}(u))$  if $u\in(0,1)^m$ and $0$ otherwise.
By \eqref{chgtvar} and the integrability assumption made on $h$, the
function $\tilde{h}$ is square integrable on $\rset^m$. Using the
definition of $\tilde{h}$ for the first equality and symmetry for the
second one, one has
\begin{multline*}
  \sum_{\bi \in \mathcal{I}_m} R_{\bi}[h,h]  =I^m \ \sum_{\bi \in \mathcal{I}_m} \int_{\mathcal{J}_{\bi}^2}\tilde{h}(u)\{\tilde{h}(u)-\tilde{h}(v)\}dudv
 = \frac{I^m}{2}  \sum_{\bi \in \mathcal{I}_m} \int_{\mathcal{J}_{\bi}^2} \{\tilde{h}(u)-\tilde{h}(v)\}^2 \ du dv \\
 = \frac{I^m}{2} \sum_{\bi \in \mathcal{I}_m} \int_{\mathcal{J}_{\bi}} \int_{\mathcal{J}_{\bi}-u}     \{\tilde{h}(u)-\tilde{h}(u+w)\}^2 \ dw du
   \leq \frac{1}{2} \int_{[0,1]^m}\int_{[-1,1]^m}(\tilde{h}(u)-\tilde{h}(u+z/I))^2 \ dudz
  \eqsp.
\end{multline*}
where we have set, for $\bi \in \mathcal{I}_m$, $\mathcal{J}_{\bi}= \prod_{k=1}^m  [(i_k-1)/I, i_k/I]$.
By continuity of the translations in $L^2(\rset^m,du)$ and Lebesgue's Theorem, one obtains that the
right-hand-side converges to $0$ as $I\to\infty$.
\end{proof}
We now proceed to the proof of Proposition \ref{prop:multVarLimGalCase}. Under A\ref{A1}, it holds that
\begin{equation}\label{minoq}
   q_{\bi}(\chi) \geq
\left( \essinf{g \cdot \leb}{ \chi /g} \right) \ \int_{\strata{\bi} }g \ d\leb=
 I^{-m} \  \essinf{g \cdot \leb}{ \chi /g} \eqsp.
\end{equation}
Hence, by Lemma~\ref{arrondis}(\ref{arrondis1}), to prove the first assertion,
it is enough to check that $\lim_{I\rightarrow +\infty}\sum_{\bi \in \{1,
  \dots, I\}^m} \frac{p_{\bi}^2 \; \sigma_{\bi}^2}{q_{\bi}(\chi)}=
\varsigma^2_\infty(\mu,\chi)$.  By definition of $R_{\bi}$ (see
Eq.~(\ref{eq:DefiResteR})),
\begin{align*}
  \frac{p_{\bi}^2 \; \sigma_{\bi}^2}{q_{\bi}(\chi)}
  &=\frac{\int_{\strata{\bi}} \marg^2(\mmarg-\smarg^2) /g \ d\leb
    -R_{\bi}[\marg,\mmarg\marg]+ R_{\bi}[\smarg\marg,\smarg\marg]}{I^m \
    \int_{\strata{\bi}}\chi \ d\leb} \eqsp,
\\\varsigma^2_\infty(\mu,\chi) &=\sum_{\bi \in \mathcal{I}_m}\frac{\int_{\strata{\bi}}
    \marg^2(\mmarg-\smarg^2) /g \ d \leb-R_{\bi}[\chi, \marg
    ^2(\mmarg-\smarg^2)/\chi]}{I^m \ \int_{\strata{\bi}} \chi \ d \leb}.\end{align*}
\[
\mbox{Therefore }\sum_{\bi \in \mathcal{I}_m} \frac{p_{\bi}^2 \;
  \sigma_{\bi}^2}{q_{\bi}(\chi)}-\varsigma^2_\infty(\mu,\chi) =\sum_{\bi
  \in\mathcal{I}_m} \frac{R_{\bi}[\chi, \marg ^2(\mmarg-\smarg^2)
  /\chi]+ R_{\bi}[\smarg\marg,\smarg\marg]-R_{\bi}[\marg,\mmarg\marg]}{I^m
  \ \int_{\strata{\bi}}\chi \ d\leb}, \] and one easily concludes with
\eqref{minoq} and Lemma \ref{lem:ControleRestes} (which applies under A\ref{A2}
and A\ref{A3}).  The second assertion is a consequence of Lemma~\ref{arrondis}(\ref{arrondis2}).
\end{proof}

\begin{proof}{\emph{ of Proposition~\ref{varasopt}}}
Since for $a,b\geq 0$, $|\sqrt{a}-\sqrt{b}|\leq \sqrt{|a-b|}$, one has
\begin{align*}
\sum_{\bi \in \mathcal{I}_m} \bigg| p_{\bi} \sigma_{\bi} &
  -\int_{\strata{\bi}} \left[\marg \sqrt{\mmarg-\smarg^2} \right] \ d\leb\bigg| \\
  &\leq \sum_{\bi \in \mathcal{I}_m} \bigg|\int_{\strata{\bi}}\marg \ d\leb\int_{\strata{\bi}}\mmarg\marg \ d\leb -\bigg(\int_{\strata{\bi}}\smarg\marg \ d\leb \bigg)^2-\bigg(\int_{\strata{\bi}}\left[\marg\sqrt{\mmarg-\smarg^2} \right]\ d\leb\bigg)^2\bigg|^{1/2}\\
  &=\sum_{\bi \in \mathcal{I}_m}\sqrt{\frac{1}{I^m}\left|-R_{\bi}[\marg,\mmarg\marg]+R_{\bi}[\smarg\marg,\smarg\marg]+R_{\bi}[\marg\sqrt{\mmarg-\smarg^2},\marg\sqrt{\mmarg-\smarg^2}]\right|}\\
  &\leq \left(\sum_{\bi \in \mathcal{I}_m}\left|-R_{\bi}[\marg,\mmarg\marg]+R_{\bi}[\smarg\marg,\smarg\marg]+R_{\bi}[\marg\sqrt{\mmarg-\smarg^2},\marg\sqrt{\mmarg-\smarg^2}]\right|\right)^{1/2}
  \eqsp.
\end{align*}
Under A\ref{A2}, $\int \marg^2(\mmarg-\smarg^2) /g \ d\leb<+\infty$, and by
Lemma \ref{lem:ControleRestes}, the right-hand-side converges to $0$ as
$I\to+\infty$. Therefore,
\begin{equation}
\label{eq:varasopt}
\lim_{I \to +\infty} \sum_{\bi \in \mathcal{I}_m} \bigg| p_{\bi}
\sigma_{\bi} -\int_{\strata{\bi}} \left[\marg \sqrt{\mmarg-\smarg^2} \right] \
d\leb\bigg| =0 \eqsp.
\end{equation}
We now write
  \begin{multline*}
    \left(\int_{\R^m} \left[\marg \sqrt{\mmarg-\smarg^2} \right] \ d\leb \right) \
    \sum_{\bi \in \mathcal{I}_m} \left|q_\bi(\chiopt{}) -
      \qopt{\bi}(\mu)
    \right| \\
    \leq \sum_{\bi \in \mathcal{I}_m} \qopt{\bi}(\mu) \left|
      \sum_{\bj \in \mathcal{I}_m} p_\bj \sigma_\bj - \int \left[\marg
        \sqrt{\mmarg-\smarg^2} \right] \ d\leb \right|
    + \sum_{\bi \in \mathcal{I}_m} \left| p_\bi \sigma_\bi -
      \int_{\strata{\bi}} \left[\marg \sqrt{\mmarg-\smarg^2} \right] \
      d\leb\right| \eqsp.
  \end{multline*}
  By Eq.\eqref{eq:varasopt}, the rhs tends to zero as $I\to +\infty$. The
  second assertion is a consequence of Lemma~\ref{arrondis}(\ref{arrondis3})
  applied with $\epsilon = \sqrt{M/I^m}$ and of Eq.~\eqref{eq:varasopt}.
\end{proof}

\subsection{Proofs of Section~\ref{sec:AdaptiveAlgorithm}}
\label{sec:Proof:Algo}
We only give the proof of Proposition~\ref{prop:DeriveeMu} and refer to
\cite{etore:fort:jourdain:moulines:2008} for the one of Corollary~\ref{coro:Derivee-Matrice}.\begin{proof}{\emph{of Proposition~\ref{prop:DeriveeMu}}}
  Let $H\in\rset^d$ be such that $|H|<|\mu|$, $e_1=\frac{\mu}{|\mu|}$,
  $a=\pscal{H}{e_1}$, $b=|H- a e_1|$ and $e_2$ be equal to $\frac{H - a
    e_1}{b}$ if $b\neq 0$ and to any vector with norm $1$ orthogonal to $e_1$
  otherwise. We complete $(e_1,e_2)$ with $(e_3,\hdots,e_d)$ to obtain an
  orthonormal basis of $\R^d$.  For $\alpha \in \R^d$, $\alpha_k =
  \pscal{\alpha}{e_k}$.
\begin{align}
\label{eq:difference-gz}
\nonumber  g_z(\mu+H)-g_z(\mu)
  &= \int_{\{\alpha, \alpha_1 \leq\frac{z-\alpha_2 b}{|\mu|+a} \} } h(\alpha) \; d\alpha - \int_{\{\alpha, \alpha_1 \leq \frac{z}{|\mu|} \}} h( \alpha ) \; d\alpha =\int_{\R^{d-1}}\int_{\frac{z}{|\mu|}}^{\frac{z-\alpha_2 b}{|\mu|+a}} h (\alpha)d\alpha_1 d\alpha_{2:d}\\
\nonumber  &=-\int_{\R^{d-1}}\int_0^1 h\left(\frac{z-\alpha_2bs}{|\mu|+as}e_1+\sum_{k=2}^d\alpha_ke_k\right)\frac{az+\alpha_2b|\mu|}{(|\mu|+as)^2}dsd\alpha_{2:d}\\
\nonumber
  &  =-\int_0^1\int_{\R^{d-1}}
  h\left(z\frac{(|\mu|+as)e_1+bse_2}{(|\mu|+as)^2+(bs)^2}+\sum_{k=3}^d\alpha_ke_k  \phantom{+\left(\alpha_2-\frac{zbs}{(|\mu|+as)^2+(bs)^2}\right)\frac{(|\mu|+as)e_2-bse_1}{|\mu|+as}}\right.\\
\nonumber
  & \left. \phantom{\sum_{k=3}^d\alpha_ke_k }+ \left(\alpha_2-\frac{zbs}{(|\mu|+as)^2+(bs)^2}\right)\frac{(|\mu|+as)e_2-bse_1}{|\mu|+as} \right)\frac{az+\alpha_2b|\mu|}{(|\mu|+as)^2}d\alpha_{2:d}ds\\
  &=-\int_0^1\int h(y)\frac{\pscal{y}{H}}{|\mu+sH|}d\lambda^{\mu+sH}_zds \eqsp,
\end{align}
where, for the last equality, we made the change of variable
$$\beta_2=\frac{\sqrt{(|\mu|+as)^2+(bs)^2}}{|\mu|+as}\alpha_2-\frac{zbs}{(|\mu|+as)\sqrt{(|\mu|+as)^2+(bs)^2}} \eqsp,
$$
used the equality $(|\mu|+as)e_1+bse_2=\mu+sH$ and remarked that
$\pscal{\mu+sH}{y}=z$ implies that $az+
\pscal{y}{e_2}b|\mu|=(|\mu|+as)\pscal{y}{H}$. Define, for $\nu\in\R^d_*$,
$\gamma(h,\nu) \eqdef \int \frac{y}{|\nu|} h(y)d\lambda^{\nu}_z$.
We deduce that
$$g_z(\mu+H)-g_z(\mu)+ \pscal{H}{\int \frac{y}{|\mu|} \; h\bigl( y\bigr) \;
  d\lebh{z}{\mu}}=\pscal{H}{\int_0^1 \{ \gamma(h,\mu)-\gamma(h,\mu+sH) \} \, ds}\eqsp.$$
Consider now the following decomposition
$\gamma(h,\nu)= \gamma\left( h \1_{\{ |\cdot| > M\}},\nu\right) + \gamma\left(h \1_{\{ |\cdot| \leq M\}},\nu\right)$.
Under assumption \eqref{eq:uniforme-integrabilite}, the first term in the
rhs is arbitrarily small as $M$ goes to infinity uniformly in $\nu$ close to $\mu$.
When $\nu \to \mu$, the measure $\1_{\{ |\cdot| \leq M\}} \lambda_z^{\nu}$ converges weakly to $\1_{\{ |\cdot| \leq M\}} \lambda_z^\mu$;
hence, the second term converges to $\gamma\left(h \1_{\{ |\cdot| \leq M\}},\mu\right)$. Therefore,
the function $\nu \mapsto \gamma(h,\nu)$ is continuous at $\mu$ and the conclusion follows easily.
\end{proof}

\begin{table}[htbp]
\centering
  \begin{tabular}{|ll|ll|c|c|c|cccc|} \hline
\multicolumn{2}{|c|}{Model} &  &  &  Price & \multicolumn{6}{c|}{Variance} \\
\hline
$\upsilon$ & $
K$ & $\nu$&  $q_i$ & - &  {\tt MC} &  {\tt AdaptStr} &  {\tt GHS}  & $\mu_{\mathrm{reg}}$  &  $\mu_{\star}$ & $\mu_{l}$  \\ \hline
0.1 & 45 &  0 & prop & 6.05 & 8.640 & - & - &0.017 & 0.016 & 0.017\\
& &   & opt & 6.05 & 8.640 & 0.004& -& 0.005 &0.004 & 0.004 \\
& &  $\nu_\star$ & prop &6.05  & 0.803 & -  & 0.014 & 0.014 & 0.008 & 0.007 \\
& &   & opt & 6.05 & 0.803  & 0.002 &- & 0.005 & 0.002  & 0.002 \\ \hline
0.5 & 45 & 0 & prop & 9.00 & 158.1  & - & - & 2.086 & 2.128 & 2.243 \\
& &  & opt & 9.00 & 158.1 & 0.352 & - & 0.362 & 0.371 & 0.390 \\
& & $\nu_\star$ & prop &9.00 & 14.95 & -  & 0.203 & 0.225 & 0.324  & 0.221 \\
 & &  & opt & 9.00 & 14.95 & 0.147 & - & 0.162 & 0.223 & 0.162 \\ \hline
0.5 & 65 & 0 & prop & 2.16 & 48.41 & - & - & 1.857  & 1.859 & 2.097\\
& & & opt & 2.16 & 48.41  & 0.093 & - & 0.096 & 0.096 &0.147 \\
& & $\nu_\star$ & prop & 2.16  & 2.32 & - &0.039 &0.046 &0.048 & 0.049 \\
& & & opt & 2.16& 2.32  & 0.020 & -  &0.024 & 0.025 & 0.026 \\ \hline
1 & 45 & 0 & prop & 14.01 & 852.0 & - & - & 52.24 & 54.51 &57.72 \\
& &  & opt &  14.01 & 852.0 & 5.39 & - & 5.48 & 5.69  & 6.02 \\
 &  & $\nu_\star$ & prop & 14.01 & 42.76 & - & 3.014 & 3.185 & 4.360 & 3.265 \\
& & & opt & 14.01 & 42.76 & 2.270 & - & 2.400 & 3.220 & 2.450 \\ \hline
1 & 65 & 0 & prop & 7.79  & 587 & -  & -  & 50.9& 50.5 & 55.8 ; \\
& & & opt & 7.79 & 587  &3.75 & -  & 3.01 & 3.08 & 3.95 \\
 &  & $\nu_\star$ & prop &7.78&  22.34 & - & 1.55 & 1.75 & 2.01 & 1.56 \\
& & & opt & 7.78& 22.34   & 0.99 & - & 1.14  & 1.31 & 1.00\\ \hline
\end{tabular}
\caption{Asian Option: Monte Carlo and stratification}
\label{tab:Asie}
\begin{center}
  \begin{tabular}{|ll|l|c| c c c c|} \hline
\multicolumn{2}{|c|}{Model} &  &   Price & \multicolumn{4}{c|}{Variance} \\
$\upsilon$ & $K$ & $\nu$ & & Latin  & Latin +Rot $\mu_{\mathrm{reg}}$ & Latin +Rot $\mu_{\star}$ & Latin +Rot $\mu^{(N)}$\\ \hline
0.1& 45 & 0 & 6.05 &0.0596 & 0.0008)  & 0.0008 & 0.0008\\
& & $\nu_\star$ & 6.05 &  0.6000 & 0.0063 & 0.0009 & 0.0003   \\ \hline
0.5 & 45 & 0 & 9.00 & 35.55 & 0.374  & 0.351  & 0.385  \\
& & $\nu_\star$ & 9.00 &  11.72 & 0.166 & 0.242 &0.137\\ \hline
0.5 & 65 & 0 &  2.16 & 27.55 &0.152 & 0.135 & 0.147\\
 & & $\nu_\star$ & 2.16 & 2.00 & 0.043 &  0.037 & 0.033\\ \hline
1 & 45 & 0 & 14.00 &  357.70 & 10.86  & 9.84  & 12.20 \\
& & $\nu_\star$ & 14.00 & 36.25 & 2.35 & 3.25  & 2.10  \\ \hline
1 & 65 & 0 & 7.78 & 339.11 & 7.94 & 7.70 & 9.14 \\
& & $\nu_\star$ & 7.78 & 19.62 & 1.49 & 1.34 & 1.25 \\
\hline \end{tabular}
\caption{Asian Option: Latin Hypercube }
\label{tab:AsieLHS}
\end{center}
\end{table}

\begin{table}[htbp]
\centering
  \begin{tabular}{|ll|ll|c|c|c|cccc|}\hline
\multicolumn{2}{|c|}{Model} &  &  &  Price & \multicolumn{6}{c|}{Variance} \\
\hline
$K$ & $B$ & $\nu$&  $q_i$ & - &  {\tt MC} &  {\tt AdaptStr} &  {\tt GHS}  & $\mu_{\mathrm{reg}}$  &  $\mu_{\star}$ & $\mu_{l}$  \\
\hline
50 & 60 & 0 & prop & 1.38 & 2.99 & - & - & 1.46 & 1.13 & 1.14 \\
& &  & opt & 1.38 & 2.99 & 0.31  & - &0.83 & 0.31 & 0.31 \\
 &  & $\nu_\star$  & prop & 1.38 & 1.34  & - & 0.50 & 1.15 & 0.49 & 0.50 \\
& & & opt &1.38 & 1.34 & 0.17 & - & 1.12 & 0.31 & 0.31 \\ \hline
50 & 80 & 0 & prop & 1.92 & 4.92 & - & - & 0.016  & 0.017 & 0.016 \\
 &  &  & opt &1.92 & 4.92 & 0.002 & - & 0.002 & 0.002 &0.002 \\
& & $\nu_\star$  & prop & 1.92 & 0.704 & - & 0.0011 &  0.0012  & 0.0013   &  0.0011 \\
& & & opt & 1.92 & 0.704 & 0.0005 & - &0.0006 & 0.0006 & 0.0005\\ \hline \end{tabular}
\caption{Barrier Option: Monte Carlo and stratification}
\label{tab:Barrier}

\begin{center}
  \begin{tabular}{|ll|l|c| cccc|} \hline
\multicolumn{2}{|c|}{Model} &  &   Price & \multicolumn{4}{c|}{Variance} \\
$K$ & $B$ & $\nu$ & & Latin  & Latin +Rot $\mu_{\mathrm{reg}}$ &  Latin +Rot $\mu_\star$  &  Latin +Rot $\mu^{(N)}$  \\ \hline
50 & 60 & 0 & 1.38 & 1.98 & (1.5074; 1.5416) & 0.97 & 0.98 \\
& & $\nu_\star $ & 1.38 & 1.26 & 1.01& 0.31 & 0.21 \\ \hline
50 & 80 & 0 &  1.92 & 0.727 & 0.002 & 0.002 & 0.002 \\
& &  $\nu_\star $ & 1.92 & 0.4501& 0.0005 & 0.0006 & 0.0004  \\
\hline \end{tabular}
\caption{Barrier Option: Latin Hypercube }
\label{tab:BarrierLHS}
\end{center}

\centering
   \begin{tabular}{|ll|ll|c|c|c|cccc|}\hline
\multicolumn{2}{|c|}{Model} &  &  &  Price & \multicolumn{6}{c|}{Variance} \\
\hline
$c$ & $K$ & $\nu$&  $q_i$ & - &  {\tt MC} &  {\tt AdaptStr} &  {\tt GHS}  & $\mu_{\mathrm{reg}}$  &  $\mu_{\star}$ & $\mu_{l}$  \\
\hline
0.1 & 45 & 0 & prop & 11.24 & 22.16 & - & - & 0.25 & 0.25 & 0.25 \\
& & & opt & 11.24 & 22.16 & 0.22 & - & 0.22; & 0.22 & 0.22 \\
& & $\nu_\star$ & prop & 11.24 & 1.39 & - & 0.26 & 0.87 & 0.26 & 0.26 \\
& & & opt &   11.24 & 1.39 & 0.21 & -  & 0.65 & 0.21 & 0.21\\ \hline
0.5 & 45 & 0 & prop & 11.56 & 51.13 & - & - & 0.37 & 0.37 & 0.37 \\
& & & opt & 11.56 & 81.13 & 0.10 & - & 0.10 & 0.10 & 0.10 \\
& & $\nu_\star$ & prop &  11.56 & 8.64 & - & 0.08 & 0.09 & 0.09 & 0.08  \\
& & & opt &  11.56 & 8.64 &0.06& - & 0.07 & 0.07 & 0.06 \\ \hline
0.9 & 45 & 0 & prop & 12.09 & 134 & - & - & 0.75 &  0.74 &  0.74 \\
& &  & opt & 12.09 & 134 &0.05 & - & 0.05 & 0.05 & 0.05 \\
& & $\nu_\star$ & prop & 12.09 & 14.46 & - & 0.022 & 0.029 & 0.024 & 0.023 \\
&  &  & opt & 12.09 & 14.46& 0.008 & - & 0.012 & 0.009 & 0.008 \\
\hline \end{tabular}
\caption{Basket  Option: Monte Carlo and stratification}
\label{tab:Basket}

\begin{center}
  \begin{tabular}{|ll|l|c|cccc|} \hline
\multicolumn{2}{|c|}{Model} &  &   Price & \multicolumn{4}{c|}{Variance} \\
$c$ & $K$ & $\nu$ & & Latin &   Latin +Rot $\mu_{\mathrm{reg}}$ &  Latin +Rot $\mu_\star$  &  Latin +Rot $\mu^{(N)}$  \\ \hline
0.1 & 45 & 0 & 11.24 & 0.08 & 0.07 & 0.03 &0.03 \\
& & $\nu_\star$ & 11.24 &  1.18 &0.92 &0.05 & 0.04 \\ \hline
0.5 & 45 & 0 & 11.56 & 4.94 & 0.02 & 0.02 & 0.02 \\
& & $\nu_\star$ &11.56 & 6.90 & 0.02 & 0.02 & 0.02 \\ \hline
0.9 & 45 & 0 & 12.09 & 13.05 & 0.007 & 0.006 & 0.007  \\
& & $\nu_\star$ &12.09 & 12.51 & 0.0038 & 0.0026 & 0.0006\\ \hline
\end{tabular}
\caption{Basket Option: Latin Hypercube }
\label{tab:BasketLHS}
\end{center}
\end{table}

\begin{table}[htbp]
\centering
   \begin{tabular}{|ll|ll|c|c|c|c|}\hline
\multicolumn{2}{|c|}{Model} &  &  &  Price & \multicolumn{3}{c|}{Variance} \\
\hline
$\xi_0$ & $K$ & $\nu$&  $q_i$ & - &  {\tt MC} &  {\tt AdaptStr} &  $\mu_{\mathrm{reg}}$    \\
\hline
0.01 & 120 & 0 & prop &  0.007  & 0.0272 &  - &  0.0234 \\
 &  &  & opt &   0.007  &0.0272 &  0.0003  & 0.0003 \\
 & 100 & 0 & prop & 5.19  & 18.06 &  - & 2.009  \\
& & & opt &  5.19  & 18.06  & 1.700 & 1.725 \\
 & 80 & 0 & prop &  22.65 & 30.01 &  - & 3.16  \\
& & & opt &  22.65  & 30.01  & 2.85 & 2.93 \\ \hline
0.04& 130 & 0 & prop & 0.024  & 0.152 &  - &  0.098  \\
 &  &  & opt & 0.024  & 0.152 &0.001   &0.001  \\
& 100 & 0 & prop & 6.42 & 46.58 & -  & 2.37 \\
& & & opt & 6.42 & 46.58 & 1.69 & 1.68 \\
 & 70 & 0 & prop & 31.74 & 88.38 & - & 4.08 \\
& & & opt & 31.74 & 88.38 & 3.58 & 3.71 \\ \hline
\end{tabular}
\caption{Asian Option in Heston model: Monte Carlo and stratification}
\label{tab:Heston}
\end{table}
\begin{table}[htbp]
\begin{center}
  \begin{tabular}{|ll|l|c| c c c|} \hline
\multicolumn{2}{|c|}{Model} &  &   Price & \multicolumn{3}{c|}{Variance} \\
$\xi_0$ & $K$ & $\nu$ &  & Latin &  Latin +Rot $\mu_{\mathrm{reg}}$ &  Latin +Rot $\mu^{(N)}$   \\ \hline
0.01  & 120 & 0 & 0.007 & 0.027  & 0.008 & 0.009 \\
 & 100 & 0 & 5.19 & 2.77  & 2.08 & 2.15   \\
 & 80 & 0 & 22.65 & 3.65 & 2.37 & 2.20   \\ \hline
0.04 & 130 & 0 & 0.024 & 0.15 & 0.03 & 0.03\\
& 100 & 0 & 6.42 & 7.58 & 2.59 & 3.22 \\
& 70 & 0 & 31.74 & 3.19 & 3.81 &3.54 \\ \hline
\end{tabular}
\caption{Asian Option in Heston model: Latin Hypercube  }
\label{tab:HestonLHS}
\end{center}
\end{table}

\end{document}